\documentclass[reqno,11pt]{amsart}
\usepackage[utf8]{inputenc}
\usepackage{graphicx}
\usepackage{amssymb}
\usepackage{color}
\usepackage{slashed}
\usepackage{enumitem}
\usepackage[mathscr]{eucal}
\usepackage[usenames,dvipsnames]{pstricks}
\usepackage{pst-grad} 
\usepackage{pst-plot} 

\numberwithin{equation}{section}
\allowdisplaybreaks[1]

\newcommand{\SetFigFont}[3]{}

\title[The Causal Variational Principle on the Sphere]{Singular Support of Minimizers of
the \\ Causal Variational Principle on the Sphere}

\author[L.\ B\"auml]{Lucia B\"auml}
\address{Department Mathematik \\ FAU Erlangen-N\"urnberg \\ D-91058 Erlangen \\ Germany}
\email{lucia.baeuml@fau.de}

\author[F.\ Finster]{Felix Finster}
\address{Fakult\"at f\"ur Mathematik \\ Universit\"at Regensburg \\ D-93040 Regensburg \\ Germany}
\email{finster@ur.de}

\author[D.\ Schiefeneder]{Daniela Schiefeneder}
\address{Institut f\"ur Mathematik \\ Universit\"at Innsbruck \\ A-6020 Innsbruck \\ Austria}
\email{Daniela.Schiefeneder@uibk.ac.at}

\author[H.\ von der Mosel]{Heiko von der Mosel \\ \\ August 2018 / August 2019}
\address{Institut f\"ur Mathematik  \\ RWTH Aachen \\ D-52062 Aachen \\ Germany}
\email{heiko@instmath.rwth-aachen.de}

\newtheorem{Def}{Definition}[section]
\newtheorem{Thm}[Def]{Theorem}
\newtheorem{Prp}[Def]{Proposition}
\newtheorem{Lemma}[Def]{Lemma}

\newtheorem{Corollary}[Def]{Corollary}

\newcommand{\Thanks}{\vspace*{.5em} \noindent \thanks}
\newcommand{\beq}{\begin{equation}}
\newcommand{\eeq}{\end{equation}}
\newcommand{\Proof}{\begin{proof}}
\newcommand{\QED}{\end{proof} \noindent}

\newcommand{\C}{\mathbb{C}}
\newcommand{\R}{\mathbb{R}}

\newcommand{\N}{\mathbb{N}}

\newcommand{\dd}{{\mathsf{d}}}
\DeclareMathOperator{\supp}{supp}
\newcommand{\F}{{\mathscr{F}}}

\reversemarginpar

\DeclareMathOperator{\Tr}{Tr}

\renewcommand{\O}{{\mathscr{O}}}
\renewcommand{\L}{{\mathcal{L}}}
\newcommand{\Sact}{{\mathcal{S}}}

\renewcommand{\H}{{\mathscr{H}}}
\newcommand{\D}{{\mathscr{D}}}

\newcommand{\I}{{\mathcal{I}}}
\newcommand{\J}{{\mathcal{J}}}
\newcommand{\K}{{\mathcal{K}}}
\newcommand{\lbra}{\langle}
\newcommand{\lket}{\rangle}

\newcommand{\bpm}{\begin{pmatrix}}
\newcommand{\epm}{\end{pmatrix}}
\newcommand{\la}{\langle}
\newcommand{\ra}{\rangle}

\begin{document}

\begin{abstract}
The support of minimizing measures of the causal variational principle on the sphere is analyzed.
It is proven that in the case~$\tau > \sqrt{3}$,
the support of every minimizing measure
is contained in a finite number of real analytic curves which intersect at a finite number of points.
In the case~$\tau > \sqrt{6}$, the support is proven to have Hausdorff dimension 
at most~$6/7$.
\end{abstract}

\maketitle
\tableofcontents
\section{Introduction and Main Results}
In the physical theory of causal fermion systems,
space-time and the structures therein are described by a minimizer
of the so-called causal action principle
(for an introduction and the physical context see the textbook~\cite{cfs} or the survey article~\cite{dice2014}).
{\em{Causal variational principles}} were introduced in~\cite{continuum}
as a mathematical generalization of the causal action principle (for a recent introduction see~\cite{jet}).
The starting point is a topological manifold~$\F$
and a non-negative function~$\L : \F \times \F \rightarrow \R^+_0$ (the {\em{Lagrangian}}) which is assumed
to be lower semi-continuous.
The causal variational principle is to minimize the {\em{action}}~$\Sact$ defined as
the double integral over the Lagrangian
\beq \label{Sact} 
\Sact (\rho) = \int_\F d\rho(x) \int_\F d\rho(y)\: \L(x,y)
\eeq
under variations of the measure~$\rho$ within the class of regular Borel measures,
keeping the total volume~$\rho(\F)$ fixed ({\em{volume constraint}}).
Given a minimizer~$\rho_{\min}$, 
the support of the measure is referred to as
\beq \label{Mdef}
\text{space-time} \qquad M:= \supp \rho_{\min} \:.
\eeq
Structures in space-time are obtained by restricting corresponding structures on~$\F$
to~$M$. For example, a notion of {\em{causality}} is obtained by restricting the Lagrangian
to~$M \times M$: Two space-time points~$x,y \in M$
are said to be timelike and space-like separated if~$\L(x,y)>0$ and~$\L(x,y)=0$, respectively.
For more details on this notion of causality, its connection to the causal structure
in Minkowski space and to general relativity we refer to~\cite[Chapter~1]{cfs}, \cite{nrstg}
and~\cite[Sections~4.9 and~5.4]{cfs}.

Typically, the support of a minimizing measure is {\em{singular}} in the sense that it
is a low-dimensional or even discrete subset of~$\F$.
This observation goes back to~\cite{support}, where it was shown under general assumptions
that the support has an empty interior (see~\cite[Section~3.3]{support}),
\[ 
\overset{\circ}{M} = \varnothing \:. \]
A general discreteness result would be interesting from the physical point of view
because it would imply that space-time is discrete on a microscopic
scale, thereby avoiding the ultraviolet divergences of quantum field theory
(for a discussion of this point see~\cite[Section~4]{rrev}).

The present paper is devoted to analyzing the singular structure of the support
of minimizing measures. In order to keep the setting as simple as possible,
we restrict attention to the {\em{causal variational principle on the sphere}}
as introduced in~\cite[Chapter~1]{continuum} and analyzed
in~\cite[Sections~2 and~5]{support}.
Thus we choose~$\F=S^2$. For a given parameter~$\tau \in [1,\infty)$,
we define the {\em{Lagrangian}}~$\L$ by
\beq \label{Lform}
\L = \max(0, \D) \in C^{0,1}(\F \times \F, \R^+_0) \:,
\eeq
where~$\D$ is the smooth function
\beq \label{DS2}
\D(x,y) = 
\frac 14\: \big(1+ \langle x,y \rangle \big) \left( 2 - \tau^2 \:(1 - \langle x,y \rangle) \right) ,
\eeq
and~$\langle x,y \rangle$ denotes the Euclidean scalar product of two unit
vectors~$x,y \in S^2 \subset \R^3$.
The resulting causal variational principle is indeed the simplest non-trivial special case of the
causal action principle (see Examples~1.5, 1.6 and~2.8 in~\cite{continuum}; compared
to the convention in these papers, for convenience we here multiplied the Lagrangian
by a prefactor~$1/(8 \tau^2)$).
The numerical study in~\cite[Section~2]{support} gives a strong indication
for the following
\begin{quote}
{\em{Conjecture}}: If~$\tau> \sqrt{2}$, the support of every minimizing measure
consists of a finite number of points,
\[ \# M < \infty \:. \]
\end{quote}
The present mathematical methods do not seem strong enough for proving this conjecture.
But we succeed in proving the following weaker results.
\begin{Thm} \label{thm1}
In the case~$\tau > \sqrt{3}$, the support of any minimizing measure
is contained in a finite number of real analytic curves which intersect at a finite number of points.
\end{Thm} \noindent
This result clearly implies that the Hausdorff dimension of the support is at most one.
Under a stronger assumption on~$\tau$, we even prove that the Hausdorff dimension is
strictly smaller than one:
\begin{Thm} \label{thm0}
In the case~$\tau > \sqrt{6}$,
the support of any minimizing measure is totally disconnected and has Hausdorff dimension at most $6/7$.
\end{Thm}

The paper is organized as follows. In Section~\ref{secgmc} we explain the
general mathematical context and give references to related work.
In Section~\ref{secprelim} we
recall the definition of the causal variational principle on the sphere
and collect a few results from~\cite{support}.
Section~\ref{secnodal} is devoted to the so-called ``nodal set method,''
a novel technique which makes it possible to show that in the case~$\tau > \sqrt{3}$,
the support of a minimizing measure lies locally on the zero set of a
quadratic polynomial in~$\R^3$ restricted to the sphere (see Theorem~\ref{thmpoly}).
This implies that the support of every minimizing measure
is contained in a finite number of real analytic curves which intersect at a finite number of points
(see Theorem~\ref{thm1}).
In Section~\ref{seclightcone} it is shown with various methods and for different values of~$\tau$
that if~$p$ is a point of~$M$, then there is at least one support point on the
boundary of the light cone centered at~$p$. 
In Section~\ref{sechausdorff} the concept of two-sided uniform accumulation points
of scaling~$\beta$ is introduced (see Definition~\ref{deftwosidedbeta}),
and Theorem~\ref{thm0} is proved under the additional assumption that no
such accumulation points exist for~$\beta < 1/6$.
Finally, in Section~\ref{secdisconnect} such accumulation points are indeed ruled out
in the case~$\tau > \sqrt{6}$ and~$\beta < 1/6$, thereby completing the proof of Theorem~\ref{thm0}.

\section{General Mathematical Context} \label{secgmc}
We now put the causal variational principle on the sphere~\eqref{Sact}
with~$\L$ according to~\eqref{Lform} and~\eqref{DS2} into the general mathematical context.
Variational principles for measures involving a nonlocal ``interaction potential'' have been studied
extensively in the literature. Most work has been done in {\em{Euclidean space}}, i.e.\ for
actions of the form
\[ \Sact(\rho) = \int_{\R^n} \int_{\R^n} W(x-y) \: d\rho(x) \:d\rho(y) \]
for different choices of the function~$W : 
\R^n \rightarrow \R \cup \{\infty\}$. The applications range from
physics (electrostatic repulsion or other pair potentials) over
material sciences (particles, nanoparticles) to biology (flocking and self-organization).
For a good overview we refer to the introductions in~\cite{carrillo1, carrillo3} and the references therein.
Typically, interaction potentials are smooth away from the origin, being repulsive towards the origin and attractive towards infinity.
Also in the mathematical literature, mainly such {\em{repulsive-attractive}} potentials
are studied. In~\cite{carrillo1, carrillo3} the minimizers are analyzed
for potentials which have a pole at the origin. 
In~\cite{carrillo1, carrillo2} potentials are analyzed which are bounded at the origin.
In~\cite{carrillo2} it is shown that the dimension of the support of a minimizer is directly related to the
strength of the repulsion for small distances. It is proven
that for spherically symmetric potentials with mild repulsion, i.e.\ if
\beq \label{mildrep}
W(x) = w\big(|x| \big) \qquad \text{and} \qquad 
\lim_{r \searrow 0} w'(r)\, r^{1-\alpha} \rightarrow -C
\eeq
with~$C>0$ and~$\alpha>2$, the support of every minimizing measure consists of finitely many points.

On the {\em{sphere}}, variational problems of the form~\eqref{Sact} have been studied
mainly for Riesz potentials. Restricting attention to counting measures whose support consists
of~$N$ points, one is led to the problem of finding optimal distributions of points
on a sphere (see~\cite{saff+kuijlaars, cohn-kumar} and the references therein).
In~\cite{bilyk} general measures on the sphere are considered, and it is analyzed 
for Riesz potentials whether the volume measure on the sphere is minimizing or maximizing.
The main difference between between all the above variational principles and the
{\em{causal variational principle}} is that the Lagrangian~$\L(x,.)$ in~\eqref{Lform} is not
smooth but only Lipschitz continuous, having a ``cusp-like'' behavior on the zero set of the function~$\D(x,.)$.
This zero set describes a great circle centered at~$x$, and the radius of this circle
is determined by the parameter~$\tau$ (for more details see Section~\ref{secprelim} below).
Indeed, the value of this parameter has a crucial effect on the structure of the minimizing measures.
Namely, as shown in~\cite{support}, at the value~$\tau=\sqrt{2}$ there is a {\em{phase transition}}
from so-called generically timelike minimizers to minimizers with singular support
(for another effect of phase transition in a model for flocking see~\cite{burchard, frank-lieb}).
In the present paper, we restrict attention to the parameter range~$\tau>\sqrt{2}$ where
the minimizing measures have singular support. Our goal is to derive upper bounds for the Hausdorff
dimension of the support. Therefore, our goal is similar to the results in~\cite{carrillo2}. 
But our potential has a rather different form for two reasons.
First, as explained above the potential~$\L(x,.)$ is only Lipschitz on a great circle centered at~$x$.
Second, the behavior of the potential near the origin is different. Indeed, the potential~$\L(x,y)$ is
smooth near the diagonal and has a local maximum at~$x=y$ (thus it is also repulsive for
small distances). It has an expansion similar to~\eqref{mildrep}, however with~~$\alpha=2$
instead of~$\alpha>2$ (see~\eqref{Dp1}, where~$\vartheta$ is the geodesic distance).
Due to these two major differences, we must employ very different techniques.

\section{Preliminaries} \label{secprelim}
We now recall the setup of the causal variational principle on the sphere
and review a few basic results from~\cite{support} of relevance for what follows.
Choosing~$\F = S^2$ and~$\tau \geq 1$,
we minimize the action~\eqref{Sact} for the Lagrangian~\eqref{Lform} and~\eqref{DS2},
under variations of the measure~$\rho$ within the class of normalized regular Borel measures.
The existence of minimizers follows immediately from abstract compactness arguments
(see~\cite[Section~1.2]{continuum}). The minimizers will in general not be unique.
In what follows, we always let~$\rho_{\min}$ be a minimizing measure
and denote its support by~$M$, \eqref{Mdef}.

It is convenient to denote the angle between two points~$x,y \in S^2$
by~$\vartheta_{xy} \in [0, \pi]$. Since the function~$\D$ in~\eqref{DS2} depends only
on this angle, the Lagrangian~\eqref{Lform} can be written as
\begin{align}
\L(x,y) &= \max \big(0, \D(\vartheta_{xy}) \big) \quad \text{with} \label{Lform3} \\
\D(\vartheta) &= \frac 14\: (1+ \cos \vartheta) \left( 2 - \tau^2 \:(1 - \cos \vartheta) \right) . \label{DS3}
\end{align}
The function~$\D$ has its maximum at~$\vartheta=0$. It has zeros at~$\vartheta$ equals~$\pi$ and
\beq \label{thetamax}
\vartheta_{\max} := \arccos\left(1-\frac{2}{\tau^2}\right) ,
\eeq 
being positive on the interval~$[0, \vartheta_{\max})$ and negative on the interval~$(\vartheta_{\max}, \pi)$.

Motivated by the causal structure as defined after~\eqref{Mdef}, we introduce the sets
\begin{align*}
\I(x) &= \{ y \in \F \text{ with } \D(x,y) > 0 \} &\quad& \text{open light cone} \\
\J(x) &= \{ y \in \F \text{ with } \D(x,y) \geq 0 \} && \text{closed light cone} \\
\K(x) &= \partial \I(x) \cap \partial \big( \F \setminus \J(x) \big) && \text{boundary of the light cone}\:.
\end{align*}
Given a normalized regular Borel measure~$\rho$ (not necessarily a minimizer),
we introduce the Hilbert space~$L^2(\F, d\rho)$ by~$(\H_\rho, \lbra .,. \lket_\rho)$
and define the operator~$\L_\rho$ by
\beq \label{Lrhodef}
\L_\rho \::\: \H_\rho \rightarrow \H_\rho \:,\: \qquad (\L_\rho \psi)(x) =
\int_\F \L(x,y)\: \psi(y)\: d\rho(y) \:.
\eeq
This operator is self-adjoint and Hilbert-Schmidt (see~\cite[Lemma~3.1]{support}).

We now state the Euler-Lagrange equations:
\begin{Prp} \label{prpEL}
Let~$\rho_{\min}$ be a minimizer of the causal variational principle on the sphere.
Then the function~$\ell$ defined by
\beq \label{ldef}
\ell(x) = \int_\F \L(x,y)\: d\rho_{\min}(y) \;\; \in C^{0,1}(\F)
\eeq
is minimal on~$M$,
\beq \label{EL1}
\ell|_M \,\equiv\, \inf_\F \ell \:.
\eeq
Moreover, the operator~$\L_{\rho_{\min}}$ is positive semi-definite, i.e.
\beq \label{EL2}
\lbra \psi, \L_{\rho_{\min}} \psi \lket_{\rho_{\min}} \geq 0 \qquad \text{for all~$\psi \in \H_{\rho_{\min}}$}\:.
\eeq
\end{Prp} \noindent
The equations~\eqref{EL1} are derived by considering first variations of the measure~$\rho_{\min}$
(see~\cite[Lemma~3.4]{support}). The positivity statement~\eqref{EL2} is
obtained by considering second variations and using that they are non-negative
(see~\cite[Lemma~3.5]{support}).

As a particular consequence of \eqref{EL2}, any finite set
of points in the support of the minimizing measure gives rise to a symmetric
positive semi-definite matrix:
\begin{Corollary}\label{corollary4}
Let $N\in\N$ and $p_0,\ldots,p_N\in M$, then the Gram matrix $L$ defined by
$$
L=\Big(\L(p_i,p_j)\Big)_{i,j=0,\ldots,N}\in\C^{(N+1)\times (N+1)}
$$
is symmetric and positive semi-definite.
\end{Corollary} \noindent
For the proof see \cite[Corollary 3.6]{support}.

In the subsequent analysis, we will make use of the specific form of the
function~$\D$ in~\eqref{DS2}. The key point is that this function is
a polynomial in~$\la x,y \ra$ of degree two. This makes it possible
to write it in terms of spherical harmonics~$Y_l^m$ of total angular momentum~$l \leq 2$.
More precisely, a direct computation yields (see~\cite[Section~5.1]{support})
\beq \label{DYrep}
\D(x,y)= 4\pi  \sum_{l=0}^{2}\nu_l \sum_{m=-l}^l Y_l^m(x)\, \overline{Y_l^m(y)} \:,
\eeq
where the coefficients~$\nu_l$ are given by (compared to the formulas in~\cite[eq.~(5.1)]{support},
we again multiplied by a factor~$1/(8 \tau^2)$)
\beq \label{nu0}
\nu_0= \frac 12 -\frac{1}{6}\:\tau^2 \:,\qquad\nu_1=\frac{1}{6} \qquad \text{and} \qquad
\nu_2=\frac{1}{30}\:\tau^2\:.
\eeq

We finally recall a few general results on the structure of minimizing measures.
If~$\tau \leq \sqrt{2}$, the minimizers are generically timelike,
which means in particular that all points in the support of~$\rho_{\min}$
are timelike or lightlike separated (for details see~\cite[Definition~3.8 and Corollary~5.1]{support}).
In the remaining case~$\tau> \sqrt{2}$, the support of every minimizing
measure is singular in the sense that it has an empty interior
(see~\cite[Theorem~3.18 and Corollary~5.1]{support}).

\section{The Nodal Set Method} \label{secnodal}
\begin{Def} An open subset~$\Omega \subset S^2$ is called {\bf{totally timelike}}
if~$\D(x,y) > 0$ for all~$x,y \in \Omega$.
\end{Def}
In what follows, we let~$\Omega$ be a totally timelike set. We denote
the restriction of the minimizing measure~$\rho_{\min}$ to~$\Omega$ by~$\rho$,
\[ \rho := \chi_\Omega\: \rho_{\min} \]
(where~$\chi$ is the characteristic function).
Extending functions on~$\supp \rho$ by zero to all of~$S^2$, the Hilbert space~$\H_\rho$
becomes a subspace of~$\H_{\rho_{\min}}$. As a consequence, the positivity result of
Proposition~\ref{prpEL} also applies to the measure~$\rho$,
\beq \label{ELrho}
\la \psi, \L_\rho \psi \ra_\rho \geq 0 \qquad \text{for all~$\psi \in \H_\rho$}\:.
\eeq
Combining~\eqref{Lform} with the fact that the set~$\Omega$ is totally timelike,
we may replace the kernel~$\L(x,y)$ in the integrand of~\eqref{Lrhodef}
by~$\D(x,y)$,
\[ (\L_\rho \psi)(x) = \int_\F \D(x,y)\: \psi(y)\: d\rho(y) \:. \]
Using the representation~\eqref{DYrep} in terms of spherical harmonics,
one sees that the image of the operator~$\L_\rho$ is a linear combination
of the vectors~$\chi_\Omega Y^m_l$. We denote this vector space by~$K$,
\beq \label{Kdef}
K := \text{span} \big\{ \chi_\Omega\, Y^m_l \:\big|\: 0 \leq l \leq 2\:,\; -l \leq m \leq l \big\} \:.
\eeq
Since the functions~$Y^m_l$ are real analytic, their restrictions to~$\Omega$ are linearly independent.
Therefore, considering~$K$ as a space of functions on~$\Omega$, this vector space
is nine-dimensional. The subtle point is that, when considered as a subspace of~$\H_\rho$,
the dimension of~$K$ might be smaller, because forming equivalence classes
of the functions~$\chi_\Omega Y^m_l$ in~$L^2(S^2, d\rho)$, the resulting vectors of~$\H_\rho$
are not necessarily linearly independent.

We let~$\mu$ be the Haar measure on~$S^2$ restricted to~$\Omega$.
Since the measure~$\mu$ is supported in all of~$\Omega$, considering~$K$ as a subspace
of~$\H_\mu$, this vector space is indeed nine-dimensional.
In analogy to~\eqref{Lrhodef} we introduce the operator
\beq \label{Lmudef}
L_\mu \::\: K \subset \H_\mu \rightarrow K \:,\qquad (L_\mu \psi)(x) = \chi_\Omega(x)
\int_{\Omega} \D(x,y)\: \psi(y)\: d\mu(y) \:.
\eeq
Clearly, the operator~$L_\mu$ maps to~$K$. Therefore, being the restriction of a symmetric operator
to an invariant subspace, the operator~$L_\mu$
is again symmetric. We now determine the signs of its eigenvalues:

\begin{Lemma} \label{lemmamultiply}
If~$\tau > \sqrt{3}$, the operator~$L_\mu$ has one negative
and eight positive eigenvalues (where we count the eigenvalues with multiplicities).
\end{Lemma}
\Proof We choose functions~$(\psi^m_l)$ in~$K$ for~$0 \leq l \leq 2$ and~$-l \leq m \leq l$ such that
\[ \la \psi^m_l, \chi_\Omega Y^{m'}_{l'} \ra_\mu = \delta^{m,m'}\, \delta_{l,l'}\:. \]
Then, using the definition of the operator~$L_\mu$ together with~\eqref{DYrep}, one sees that
\[ \la \psi^m_l, L_\mu \,\psi^{m'}_{l'} \ra_\mu
= 4 \pi\, \nu_\ell\: \delta^{m,m'}\, \delta_{l,l'} \:. \]
The explicit formulas~\eqref{nu0} show that in the case~$\tau>\sqrt{3}$,
the parameter~$\nu_0$ is negative, whereas~$\nu_1$ and~$\nu_2$ are positive.
Therefore, the bilinear form~$\la . , L_\mu . \ra : K \times K \rightarrow \C$
has signature~$(8,1)$. As a consequence, the operator~$L_\mu$ has one positive and
eight negative eigenvalues.
\QED

This lemma shows in particular that~$L_\mu$ is invertible. In the next lemma, we show that
the operator~$L_\mu^{-1}$ is positive semi-definite on the image of the operator~$L_\rho$
defined in analogy to~\eqref{Lmudef} by
\[ L_\rho \::\: \H_\rho \rightarrow K \:,\qquad (L_\rho \psi)(x) = \chi_\Omega(x)
\int_{\Omega} \D(x,y)\: \psi(y)\: d\rho(y) \:. \]
\begin{Lemma} \label{lemmapositive}
For any~$\psi \in \H_\rho$,
\[ \la L_\rho \psi, L_\mu^{-1} \, L_\rho \psi \ra_\mu \geq 0 \:. \]
\end{Lemma}
\Proof The positivity property~\eqref{ELrho} implies that
\[ 0 \leq \la \psi, \L_\rho \psi \ra_\rho = \la \psi, L_\rho \psi \ra_\rho \]
(note that here~$L_\rho \psi \in K$, whereas~$\L_\rho \psi$ is the
same function considered as a vector in~$\H_\rho$). Using that~$L_\mu$ is invertible,
it follows that
\begin{align*}
0 &\leq \la \psi, L_\mu L_\mu^{-1} L_\rho \psi \ra_\rho = 
\int_\Omega d\rho(x)\; \overline{\psi(x)} \:\big( L_\mu\, L_\mu^{-1} L_\rho \psi \big)(x) \\
&= \int_\Omega d\rho(x) \int_\Omega d\mu(y) \; \overline{\psi(x)} \:\D(x,y)\: \big( L_\mu^{-1} L_\rho \psi \big)(y) \\
&= \int_\Omega d\rho(x) \int_\Omega d\mu(y) \; \overline{\D(y,x)\:\psi(x)} \: \big( L_\mu^{-1} L_\rho \psi \big)(y) \\
&= \int_\Omega d\mu(y) \; \overline{L_\rho \psi(y)}\: \big( L_\mu^{-1} L_\rho \psi \big)(y)
= \la L_\rho \psi, L_\mu^{-1} \, L_\rho \psi \ra_\mu \:.
\end{align*}
This gives the result.
\QED

According to Lemma~\ref{lemmamultiply}, the operator~$L^{-1}_\mu$ has one negative eigenvalue.
Combining this fact with the result of Lemma~\ref{lemmapositive}, we conclude that
the operator~$L_\rho$ cannot be surjective. Thus there is a non-trivial vector~$v \in K$ with
\[ \la v, L_\rho \psi \ra_\mu = 0 \qquad \text{for all~$\psi \in \H_\rho$}\:. \]
As a consequence,
\begin{align*}
0 &= \la v, L_\rho \psi \ra_\mu = \int_\Omega d\mu(x) \int_\Omega d\rho(y)\;
\overline{v(x)}\: \D(x,y)\: \psi(y) \\
&= \int_\Omega d\mu(x) \int_\Omega d\rho(y)\;\overline{\D(y,x)\, v(x)}\: \: \psi(y)
= \int_\Omega d\rho(y)\;\overline{(L_\mu v)(y)}\: \psi(y)\:.
\end{align*}
Hence the non-trivial vector~$u := L_\mu v \in K$ has the remarkable property that
\[ \int_\Omega \overline{u(y)}\: \psi(y)\: d\rho(y) = 0 \qquad \text{for all~$\psi \in \H_\rho$}\:. \]
This in turn implies that the support of~$\rho$ lies on the nodal set of~$u$.
We have thus proved the following result:
\begin{Thm} \label{thmpoly}
Let~$\rho_{\min}$ be a minimizer of the causal variational principle on~$S^2$ for~$\tau > \sqrt{3}$.
Then for every totally timelike subset~$\Omega \subset S^2$
there is a non-trivial function~$u \in K$ (with~$K$ according to~\eqref{Kdef}) such that
\[ \supp \rho_{\min} \;\cap\; \Omega \subset u^{-1}(0) \cap \Omega \:. \] 
\end{Thm}

Next we make use of the fact that the space~$K$ as defined in~\eqref{Kdef} can be described without
referring to spherical harmonics as the restriction of all polynomials
in the three variables~$x$, $y$, $z$ of degree at most two to the set~$\Omega \subset S^2 \subset \R^3$.
Hence Theorem~\ref{thmpoly} tells us that the support of the minimizing measure~$\rho_{\min}$
lies inside~$\Omega$
on the intersection of the roots of two linearly independent quadratic polynomials in~$\R^3$: the function~$u$
and the function~$x^2+y^2+z^2-1$ describing the unit sphere.
This intersection can be analyzed explicitly by going through different cases
for the coefficients of the polynomial~$u$ (this method is outlined in the master's thesis
of one of the authors~\cite{prechtl}). Here we use more abstract methods from real algebraic geometry
which we learned from Tobias Kaiser:

\Proof[Proof of Theorem~\ref{thm1}] Given a totally timelike subset~$\Omega \subset S^2$, we let~$u$ be the
quadratic polynomial of Theorem~\ref{thmpoly} and denote its zero set by~$A := u^{-1}(0) \subset \R^3$.
Obviously, $\dim(A \cap S^2) \leq 2$. 
Moreover, we let~$f=x^2+y^2+z^2-1$ be the polynomial describing the sphere. 
We first prove that the following statements are equivalent:
\begin{itemize}
\item[(i)] $A=S^2$
\item[(ii)] $\dim (A \cap S^2) = 2$
\item[(iii)] There is~$\lambda \in \R \setminus \{0\}$ with~$u = \lambda f$.
\end{itemize}
Indeed, the implications (i)$\Rightarrow$(ii) and (iii)$\Rightarrow$(i) are obvious.
In order to prove the implication (ii)$\Rightarrow$(iii), we make use of the
fact that~$S^2$ is an irreducible algebraic variety and that~$A \setminus S^2$ is
an algebraic subset of~$S^2$ of codimension zero. As a consequence, $A \setminus S^2 = S^2$,
so that~$S^2 \subset A$. Hence the polynomial~$u$ vanishes on the zero set of~$f$,
i.e.\ $u \in {\mathcal{I}}({\mathcal{Z}}(f))$. Applying~\cite[Theorem~4.5.1, (v)$\Rightarrow$(ii)]{bochnak+coste+roy},
it follows that~$u \in (f)$. Since~$u$ is a quadratic polynomial, we obtain the result.

Since~$u$ and~$f$ are linearly independent, we know that~(iii) does not hold.
Hence~(ii) does not hold as well, implying that~$\dim (A \cap S^2) < 2$.
According to~\cite[Proposition~9.1.8]{bochnak+coste+roy}, there exists a
finite stratification. Using furthermore the property of semi-algebraic curves
in~\cite[Proposition 8.1.13]{bochnak+coste+roy}, it follows that~$A \cap S^2$
is a disjoint union of a finite number of points and intersection-free analytic line segments
(recall that a line segment~$C \subset \R^n$ is called intersection-free if there are~$a < c < d < b$
as well as an analytic function~$\gamma : (a,b) \rightarrow \R^n$ such that the restriction~$\gamma|_{(c,d)}$
is injective and~$\gamma((c,d)) = C$).
This concludes the proof.
\QED

\section{Support Points on the Light Cone Centered at Support Points} \label{seclightcone}
The remaining task is to show that, under the assumptions of Theorem~\ref{thm0},
the Hausdorff dimension of~$M$ is at most~$6/7$. In preparation, we
now consider the following question: Suppose that a point~$p \in S^2$ is in~$M= \supp \rho_{\min}$.
Are there necessarily points on the boundary of the light cone~$\K(p)$ which are also in~$M$?
In the case~$\tau \leq \sqrt{2}$, every minimizing measure is generically timelike
(as mentioned after~\eqref{nu0}), implying that the answer to the above question is no.
With this in mind, we may restrict attention to the parameter range~$\tau > \sqrt{2}$.
We begin with an affirmative answer to the above question under the
stronger assumption~$\tau > \sqrt{6}$:
\begin{Lemma} \label{lemma13} Assume that~$\tau > \sqrt{6}$.
If~$p \in M$, then the set~$\K(p) \cap M$ is non-empty.
\end{Lemma}
\Proof
Assume conversely that~$\K(p) \cap M = \varnothing$.
We denote the Cartesian coordinates in~$\R^3$ by~$(x,y,z)$ and also use standard
polar coordinates~$(\vartheta, \varphi)$ with~$z=\cos \vartheta$.
By a rotation of our coordinate system we can arrange that~$p$ is the point~$(1,0,0)$
on the equator.

We rotate the point~$p$ about the north pole by considering the curve
\[ p_s = (\cos s, \sin s, 0) \:. \]
Using our assumption~$\K(p) \cap M = \varnothing$ together with the fact that~$M$ is
a closed set (by definition of the support of a measure),
there is~$s_0>0$ such that for all~$s \in (-s_0, s_0)$, the set~$\K(p_s) \cap M$ is empty.
Hence for all~$s \in (-s_0, s_0)$,
\[ \ell(p_s) = \int_{\I(p)} \L\big(p_s, q \big)\: d\rho_{\min}(q) 
= \int_{\I(p)} \D\big(\vartheta_{p_s, q} \big)\: d\rho_{\min}(q) \:. \]
We now differentiate twice with respect to~$s$. Since~$\ell$ is minimal on~$M$~\eqref{EL1}, it follows that
\begin{align}
0 &\leq \frac{d^2 \ell(p_s)}{ds^2} \bigg|_{s=0} 
= \int_{\I(p)} \frac{d^2 \D\big(\vartheta_{p_s, q} \big)}{ds^2} \bigg|_{s=0} \: d\rho_{\min}(q) \notag \\
& = \int_{\I(p)} \D''\big(\vartheta_{p, q} \big)\, \left( \frac{d\vartheta_{p_s, q}}{ds} \bigg|_{s=0} \right)^2
d\rho_{\min}(q) \notag \\
&\quad+ \int_{\I(p)} \D'\big(\vartheta_{p, q} \big)\, \frac{d^2\vartheta_{p_s, q}}{ds^2} \bigg|_{s=0}\: d\rho_{\min}(q)
\:. \label{negative}
\end{align}

Differentiating~\eqref{DS3} gives
\begin{align}
\D'(\vartheta) &= -\frac 12 \: \big( 1+ \tau^2 \cos \vartheta \big) \: \sin \vartheta \label{Dp1} \\
\D''(\vartheta) &= -\frac 12 \: \big(  \cos \vartheta - \tau^2 + 2 \tau^2 \cos^2 \vartheta \big) \:. \label{Dpp}
\end{align}
By direct computation one verifies that in the considered range~$\tau>\sqrt{6}$
and for all~$\vartheta \in [0, \vartheta_{\max}]$ with~$\vartheta_{\max}$ according to~\eqref{thetamax},
the functions~$\D'(\vartheta)$ and~$\D''(\vartheta)$ are both strictly negative.

It remains to compute~$\vartheta_{p_s, q}$ and its derivatives. We parametrize~$q$ in polar coordinates as
\[ q = \big( \sin \vartheta\, \sin \varphi, \;\sin \vartheta\, \cos \varphi, \;\cos \vartheta \big) \:. \]
Then, using the sum rules,
\[ \cos \vartheta_{p_s, q} = \cos s\: \sin \vartheta\: \sin \varphi + \sin s\: \sin \vartheta\: \cos \varphi
= \sin \vartheta\, \sin (\varphi+ s) \:, \]
and thus
\begin{align*}
\frac{d\vartheta_{p_s, q}}{ds} &= \frac{d}{ds} \arccos \big(\sin \vartheta\, \sin (\varphi+ s) \big) 
= -\frac{\sin \vartheta\: \cos(\varphi+s)}{\big( 1-\sin^2 \vartheta\: \sin^2(\varphi+s) \big)^{1/2}} \\
\frac{d^2 \vartheta_{p_s, q}}{ds^2} &= 
\frac{\cos^2 \vartheta\: \sin \vartheta\: \sin(\varphi + s)}
{\big(1-\sin^2 \vartheta\: \sin^2( \varphi + s) \big)^{3/2}} > 0 \:.
\end{align*}
Combining these inequalities with the fact that~$\D'(\vartheta_{p, q}) \leq 0$ and~$\D''(\vartheta_{p, q})<0$,
we find that the expression~\eqref{negative} is less or equal to zero.
We conclude that~\eqref{negative} must vanish. This implies that~$\D'(\vartheta_{p, q})=0$
for all~$q \in M$. Hence~$\vartheta_{p, q}$ must be zero for all~$q \in M \cap \I(p)$.
In other words, the set~$M \cap \I(p)$ consists only of one point~$\{p\}$.
But if this is the case, the first integral in~\eqref{negative} is strictly negative, because
\[ \frac{d\vartheta_{p_s, p}}{ds} = \frac{d s}{ds}= 1\:. \]
This is a contradiction.
\QED

We now give two improvements of this lemma, which are valid for a larger range
of the parameter~$\tau$. These results will not be needed for the proof of Theorem~\ref{thm0}.
We state and prove them nevertheless, because the methods are of independent interest
and may be useful in the future for improving the results of the present paper.
\begin{Lemma} \label{lemma2} Assume that~$\tau > 2$.
If~$p \in M$, then the set~$\K(p) \cap M$ is non-empty.
\end{Lemma}
\Proof The method of the proof of Lemma~\ref{lemma13} no longer works,
because for~$\tau \in (2, \sqrt{6})$, the function~$\D''(\vartheta)$, \eqref{Dpp}, changes sign
on the interval~$[0, \vartheta_{\max}]$. Instead, we calculate the Laplacian of~$\D$,
\[ \Delta_p \D(p,q) = \frac{1}{\sin \vartheta} \frac{\partial}{\partial \vartheta}
\left( \sin \vartheta\: \frac{\partial \D(\vartheta)}{\partial \vartheta} \right) 
= -\frac{1}{2} \:\big( 2 \cos \vartheta - \tau^2 + 3 \tau^2 \cos^2 \vartheta \big) \:. \]
By direct computation one verifies that in the considered range~$\tau>2$,
the Laplacian of~$\D(p,q)$ is strictly negative for all points with timelike separation,
\beq \label{lapneg}
\Delta_p \D(p,q) < 0 \qquad \text{for all~$q \in \I(p)$}\:.
\eeq

Assume conversely that~$\K(p) \cap M= \varnothing$.
We introduce the measure~$\hat{\rho}$ by
\beq \label{hatrhodef}
\hat{\rho}(q) = \chi_{\I(p)}\, \rho_{\min}(q)
\eeq
and set
\beq \label{hatddef}
\hat{\dd}(r)=\int_{S^2}\D(r,q) \:d\hat{\rho}(q)\:.
\eeq
Then the functions~$\ell$ and~$\hat{\dd}$ coincide in a small neighborhood of~$p$. Hence
\[ \Delta_p \ell(p) = \Delta_p \hat{\dd}(p) = 
\int_{\I(p)} \Delta_p \D(p,q) \:d\hat{\rho}(q) \overset{\eqref{lapneg}}{<} 0 \:, \]
in contradiction to the minimality of~$\ell$ at~$p$
as obtained from the Euler-Lagrange equations~\eqref{EL1}.
\QED

We finally consider the general range~$\tau > \sqrt{2}$.
Here we use yet another method which relies on the representation of~$\D$
in terms of spherical harmonics~\eqref{DYrep}. This method makes it necessary to assume
that~$p$ is an {\em{accumulation point}} of support points.
\begin{Lemma} \label{lemma1} Assume that~$\tau > \sqrt{2}$.
If~$p$ is an accumulation point of~$M$, then the set~$\K(p) \cap M$ is non-empty.
\end{Lemma}
\Proof Assume conversely that~$\K(p) \cap M = \varnothing$.
We again introduce the measure~$\hat{\rho}$ and the function~$\hat{\dd}$
by~\eqref{hatrhodef} and~\eqref{hatddef}.
Then the functions~$\ell$ and~$\hat{\dd}$ coincide in a small neighborhood of~$p$.
We again denote the Cartesian coordinates in~$\R^3$ by~$(x,y,z)$.
By a rotation of our coordinate system we arrange that~$p$ is the north pole~$z=1$.

Using the representation~\eqref{DYrep} and~\eqref{nu0}, we obtain
\beq \label{dpoly}
\hat{\dd}(q) = 4 \pi \sum_{l=0}^2 \nu_l \sum_{m=-l}^l Y_l^m(q)\: \overline{Y_l^m(\hat{\rho})} \:,
\eeq
where~$\overline{Y_l^m(\hat{\rho})} := \int_{S^2} \overline{Y_l^m(q)} \:d\hat{\rho}(q)$.
We now parametrize the sphere near the north pole by the coordinates~$x, y$ and set~$z=\sqrt{1-x^2-y^2}$.
Being a polynomial of degree two, we can write~\eqref{dpoly} as
\[ 
\hat{\dd}(x,y) =  a_{1} \,x^2 + a_{2} \,y^2+b_{1} \,xy + b_{2} \,xz+b_{3} \,yz + c_1\,x+c_2\,y+c_3\,z+c_0 \]
with real coefficients~$a_i$, $b_i$ and~$c_i$.
Since the function~$Y_1^0 = \sqrt{\frac{3}{4 \pi}}\: z$ is positive on the support of~$\hat{\rho}$
(note that our assumption~$\tau > \sqrt{2}$ implies that~$\vartheta_{\max} < \frac{\pi}{2}$, so that~$z>0$),
we know that~$Y_1^0(\hat{\rho}) > 0$. Using in addition that the coefficient~$\nu_1$ is positive
(see~\eqref{nu0}), we conclude that~$c_3>0$.
Next, since~$\ell$ and therefore also~$\hat{\dd}$ are minimal at the north pole, the coefficients satisfy the relations
\[ c_1=-b_2 \qquad \text{and} \qquad  c_2=-b_3. \]
Furthermore, the minimality of the function~$\hat{\dd}$ implies that its Hessian is positive semi-definite, i.e.
\beq \label{semi}
\Tr \begin{pmatrix} 2a_1-c_3 & b_1 \\ b_1 & 2a_2-c_3 \end{pmatrix} \geq 0
\qquad \text{and} \qquad \det \begin{pmatrix} 2a_1-c_3 & b_1 \\ b_1 & 2a_2-c_3 \end{pmatrix} \geq 0 \:.
\eeq
Moreover, this Hessian cannot be positive definite, because otherwise~$p$ would necessarily be
an isolated point of~$M$. We conclude that the above determinant must vanish,
\beq \label{cond:det}
(2a_1-c_3)\cdot (2a_2-c_3)=b_1^2 \:.
\eeq
Combining this relation with the first inequality in~\eqref{semi}, we obtain the inequalities 
\[ 
2a_1-c_3 \geq 0 \qquad \text{and} \qquad 2a_2-c_3 \geq 0 \:. \]

We first consider the case $2a_1-c_3=0$. Then~\eqref{cond:det} implies that $b_1=0$. We
thus obtain the Taylor expansion
\[ \hat{\dd}(s,0)=c_3+c_0-\frac{b_2}{2}s^3-\frac{c_3}{8}s^4+ \O\big( s^5 \big).  \]
But since $c_3>0$, this contradicts the minimality of~$\hat{\dd}$ at~$s=0$.

In the remaining case $2a_1-c_3>0$, the relation~\eqref{cond:det} implies that
\[ a_2=\frac{b_1^2+c_3(2a_1-c_3)}{2(2a_1-c_3)} \:. \]
As a consequence, the vector
\[ u := \begin{pmatrix} -b_1  \\   2a_1-c_3 \end{pmatrix} \in \R^2 \]
lies in the kernel of the Hessian. Moreover, a short computation yields
\begin{align*}
\hat{\dd}(s u) =& c_0 +c_3+\frac{1}{2} \:\Big(b_1b_2-b_3(2a_1-c_3)\Big)\,\Big(b_1^2+2a_1-c_3\Big)^2 \, s^3\\
&-\frac{1}{8}\,c_3\,\Big( b_1^2+(c_3-2a_1)^2 \Big)^2 \, s^4 + \O\big( s^5 \big) \:,
\end{align*}
valid for all $s\in(-\varepsilon,\varepsilon)$, where $\varepsilon>0$ is chosen sufficiently small.
But since both $c_3>0$ and $ b_1^2+(c_3-2a_1)^2 >0$, this contradicts the minimality of~$\hat{\dd}$
at~$s=0$. This concludes the proof.
\QED

\section{Uniform Two-Sided Accumulation Points and Hausdorff Dimension} \label{sechausdorff}
In what follows, we let~$p$ be an accumulation point in~$M$.
According to Theorem~\ref{thm1}, there is a real analytic curve~$\gamma$
through~$p$ such that~$p$ is an accumulation point of support points on~$\gamma$.
For convenience, we parametrize~$\gamma : (-\delta, \delta) \rightarrow S^2$ by arc length,
\beq \label{gammadef}
\gamma : (-\delta, \delta) \rightarrow S^2 \qquad \text{with} \qquad \|\dot{\gamma}\| \equiv 1
\eeq
(where~$\|.\|$ denotes the Euclidean norm on~$\R^3 \supset S^2$)
and arrange that~$\gamma(0)=p$.

\begin{Def} \label{deftwosidedbeta}
The point~$p$ is a {\bf{uniform two-sided accumulation point}} on~$\gamma$
of scaling~$\beta>0$ if there exists $\varepsilon_0=\varepsilon_0(p)\in
(0,\delta)$ such that for  all $\varepsilon \in (0,\varepsilon_0]$, 
there are support
points~$\gamma(t_\pm) \in M$ with parameters in the range
\beq \label{uniformbeta}
-\varepsilon < t_- < -\varepsilon^{1+\beta} \qquad \text{and} \qquad
\varepsilon^{1+\beta} < t_+ < \varepsilon \:.
\eeq
\end{Def}

\begin{Prp} \label{prphausdorffbeta}
Assume that the set~$M:= \supp \rho_{\min}$ does not have uniform two-sided accumulation
points of scaling~$\beta>0$. Then its Hausdorff dimension is at most~$1/(1+\beta)$.
\end{Prp}
\Proof Given~$\alpha > 1/(1+\beta)$, our goal is to show that the Hausdorff measure
of~$M$ of dimension~$\alpha$ vanishes.
In view of Theorem~\ref{thm1}, it suffices to consider support points
on a finite number of smooth curves.
We denote one of these curves by~$\gamma : [0, \delta] \rightarrow S^2$.
By subdividing the smooth curves into smaller pieces (and therefore also decreasing~$\delta$),
we can arrange that there is a constant~$C$ such that the geodesic distance of the image points
is estimated from above by the distance of the parameter values, i.e.\
\[ \text{dist} \big( \gamma(\tau), \gamma(\tau') \big) \leq C\, \big| \tau-\tau' \big| \:. \]
By reparametrizing we can arrange that~$\delta=1$ (which changes the constant~$C$
in the last inequality to~$C \delta$). Moreover, since we have a finite number
of curves, it suffices to consider one of them. Thus the remaining task is to prove that the
set~$N := \gamma^{-1}(M) \subset [0,1]$ has vanishing Hausdorff measure,
\[ {\mathcal{H}}^\alpha(N) = 0 \:. \]

Let~$\sigma>0$. The assumption that~$\gamma$ has no uniform two-sided accumulation points
of scaling~$\beta>0$ means that for every~$t \in N$
there is a sequence~$(\varepsilon_\ell)_{\ell \in \N}$ with~$\varepsilon_\ell \searrow 0$
such that for every~$\ell$, there is no pair of points~$t_\pm \in N$ with
\beq \label{nopointsbeta}
-\varepsilon_\ell < t_- - t < -\varepsilon_\ell^{1+\beta} \qquad \text{and} \qquad
\varepsilon_\ell^{1+\beta} < t_+ - t < \varepsilon_\ell \:.
\eeq
We choose~$\ell$ such that~$\varepsilon_\ell < \sigma$.
We introduce the intervals~$I := (t-2\varepsilon_\ell^{1+\beta}, t+2\varepsilon_\ell^{1+\beta})$,
\[ K_- := \big(t-\varepsilon_\ell, t-\varepsilon_\ell^{1+\beta} \big) \qquad \text{and} \qquad
K_+ := \big( t+\varepsilon_\ell^{1+\beta}, t+\varepsilon_\ell \big) \:. \]
According to~\eqref{nopointsbeta}, at least one of the intervals~$K_-$ or~$K_+$ does not intersect~$N$.
Thus we can choose~$s \in \{\pm\}$ such that$~K_s \cap N = \varnothing$ and introduce the open interval
\[ J := K_s \cup I\:. \]

Varying~$t \in N$, we obtain an open covering~$(J(t))_{t \in N}$ of~$N$. Since~$N$ is compact,
there is a finite subcovering. Labeling this finite covering by~$k=1,\ldots, n$,
we have corresponding points~$t_k$, parameters~$\varepsilon_k<\sigma$ 
as well as
intervals~$I_k$ and~$J_k$. We can clearly arrange that the covering~$(J_k)_{k=1,\ldots, n}$ is
minimal in the sense that it is no longer a covering if any of the elements is taken out. Then each point
of the interval~$[-\sigma,1+\sigma]$
lies in at most two intervals of the covering~$(J_k)_{k=1,\ldots, n}$,
and thus
\beq \label{sumbound}
\sum_{k=1}^n \varepsilon_k \leq \sum_{k=1}^n \big| J_k \big| < 2\,(1+ 2 \sigma) \:.
\eeq

Since the $J_k$ cover~$[0,1]$ and the sets~$J_k \setminus I_k$ do not intersect~$N$,
the sets~$I_k$ cover~$N$. Noting that the~$I_k$ are intervals of length~$4 \varepsilon_k^{1+\beta}$,
we can estimate the Hausdorff measure of~$N$ by
\[ {\mathcal{H}}^\alpha(N) \leq \sum_{k=1}^n \Big( 4 \varepsilon_k^{1+\beta} \Big)^\alpha 
= 4^\alpha \sum_{k=1}^n \varepsilon_k^{\alpha \,(1+\beta) - 1}\: \varepsilon_k
\leq 4^\alpha \,\sigma^{\alpha \,(1+\beta) - 1} \sum_{k=1}^n \varepsilon_k \:, \]
where in the last step we used that the exponent~$\alpha \,(1+\beta) - 1$ is strictly positive.
Applying the bound~\eqref{sumbound} and taking the limit~$\sigma \searrow 0$
gives the result.
\QED

\section{Ruling out Uniform Two-Sided Accumulation Points} \label{secdisconnect}
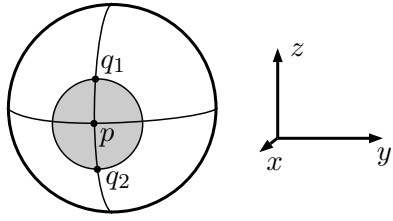
\begin{figure}
{
\begin{pspicture}(0,-1.407915)(5.0864267,1.407915)
\definecolor{colour0}{rgb}{0.8,0.8,0.8}
\pscircle[linecolor=black, linewidth=0.02, fillstyle=solid,fillcolor=colour0, dimen=outer](1.2114265,-0.21011704){0.61}
\rput[bl](1.2414266,-0.51011705){\normalsize{$p$}}
\pscircle[linecolor=black, linewidth=0.04, dimen=outer](1.4064265,0.0){1.4}
\psbezier[linecolor=black, linewidth=0.02](0.0064265444,0.0)(0.16142654,-0.13011703)(0.55142653,-0.20011704)(1.1514266,-0.2001170349121094)(1.7514266,-0.20011704)(2.6364264,-0.115117036)(2.8064265,0.0)
\psbezier[linecolor=black, linewidth=0.02](1.4064265,1.3998829)(1.2334002,1.2717129)(1.1587591,0.3630548)(1.1698854,-0.18194040728946675)(1.1810117,-0.7269356)(1.2420709,-1.2683301)(1.4064265,-1.400117)
\rput[bl](1.2514266,0.45988297){\normalsize{$q_1$}}
\rput[bl](1.3164265,-1.080117){\normalsize{$q_2$}}
\psline[linecolor=black, linewidth=0.04, arrowsize=0.05291667cm 2.0,arrowlength=1.4,arrowinset=0.0]{->}(3.6064265,-0.40011704)(5.0064263,-0.40011704)
\psline[linecolor=black, linewidth=0.04, arrowsize=0.05291667cm 2.0,arrowlength=1.4,arrowinset=0.0]{->}(3.6064265,-0.40011704)(3.6064265,0.79988295)
\psline[linecolor=black, linewidth=0.04, arrowsize=0.05291667cm 2.0,arrowlength=1.4,arrowinset=0.0]{->}(3.6014266,-0.40011704)(3.3614266,-0.58511704)
\rput[bl](3.7714264,0.694883){\normalsize{$z$}}
\rput[bl](3.4564266,-0.82011706){\normalsize{$x$}}
\rput[bl](4.9264264,-0.76511705){\normalsize{$y$}}
\pscircle[linecolor=black, linewidth=0.02, fillstyle=solid,fillcolor=black, dimen=outer](1.1639266,-0.20261703){0.0475}
\pscircle[linecolor=black, linewidth=0.02, fillstyle=solid,fillcolor=black, dimen=outer](1.1839266,0.38738295){0.0475}
\pscircle[linecolor=black, linewidth=0.02, fillstyle=solid,fillcolor=black, dimen=outer](1.2089266,-0.81261706){0.0475}
\end{pspicture}
}
\caption{The light cone~$\K(p)$ and the geometry near~$q_1$.}
\label{figsphere}
\end{figure}

In this section we complete the proof of Theorem~\ref{thm0}.
Thus we assume that~$\tau > \sqrt{6}$.
In view of Proposition~\ref{prphausdorffbeta}, we must rule out the existence
of uniform two-sided accumulation points of scaling~$\beta$ with~$\beta$ in the range
\beq \label{betacond}
0 < \beta < \frac{1}{6}\:.
\eeq
We proceed indirectly and assume that the support~$M$ of a
minimizing measure~$\rho_{\min}$ does have such a uniform two-sided accumulation point~$p$.
Then, according to Definition~\ref{deftwosidedbeta}, there is a smooth curve~$\gamma$,
parametrized by arc length and with~$\gamma(0)=p$,
such that~$p$ is an accumulation point of support points on~$\gamma$,
satisfying the uniformity condition~\eqref{uniformbeta}
for all~$\varepsilon \in (0, \varepsilon_0]$.
By a rotation in~$\R^3$ we can arrange that~$\gamma(0)= p = \big(1,0,0)$,
and that~$\dot{\gamma}(0)=(0,1,0)$ (thus the curve~$\gamma$ is tangential to the equator in the point~$p$;
see Figure~\ref{figsphere}).

\subsection{Proof that~$M \cap \K(p) \subset \{q_1, q_2\}$}
\begin{Prp} \label{lemmaq1q2}
Under the above assumptions,
\[ 
M \cap \K(p) \subset \{ q_1, q_2 \}
\qquad \text{with} \qquad
q_{1,2} = \begin{pmatrix} \cos \vartheta_{\max} \\ 0 \\ \pm \sin \vartheta_{\max} \end{pmatrix} \:. \]
\end{Prp}
\Proof Assume conversely that
there is a point~$q \in M\cap \K(p)$ with
\[ 
\langle q-p, \dot{\gamma}(0)\rangle \not=0 \]
(where~$\la .,. \ra$ denotes again the Euclidean scalar product on~$\R^3$).
Using the explicit form of $\D$ in \eqref{DS2}, its derivative is computed by
\begin{align}\label{b1}
b := \frac{d}{dt}\D \big( \gamma(t),q \big)\Big|_{t=0} & = -\D'(\vartheta_{\max})\: \frac{\langle q,\dot{\gamma}(0)\rangle}{\sqrt{1-\langle p,q\rangle^2}}
\notag\\
& = \frac{(\tau^2-1)^{3/2}}{\tau^2}\:\frac{\langle q-p,\dot{\gamma}(0)\rangle}{\sqrt{1-\langle p,q\rangle^2}}
\:\neq\: 0 \:.
\end{align}
Since $\D(p,q)=0$, we infer from \eqref{b1} that the function $t\mapsto \D(\gamma(t),q)$ 
changes sign at $t=0$. By changing the orientation of the parametrization
of $\gamma$ and decreasing~$\delta$ if necessary, we can arrange that~$\D(\gamma(t),q)>0$
and $\D(\gamma(-t),q)<0$ for all~$t \in (0, \delta)$.

Since~$p$ is a uniform two-sided accumulation point on~$\gamma$ of scaling~$\beta$,
(see Definition~\ref{deftwosidedbeta}), for any sufficiently
small~$\varepsilon<\delta$, there are parameter values~$t_\pm$ in the range~\eqref{uniformbeta}
such that the points~$\gamma(t_\pm)$ are in~$M$.
We now consider the Gram matrix~$L$ (see Corollary~\ref{corollary4} for~$N=3$) for the four points
\[ p_0 =\gamma(t_+)\:,\qquad p_1 =p \:,\qquad p_2 =\gamma(t_-)
\qquad \text{and} \qquad p_3 =q \:. \]
Then the open light cone centered at~$p_3$ contains only~$p_0$ (but not~$p_1$ and~$p_2$),
and from~\eqref{b1} we conclude that
\[ \L(p_0,p_3)=b\,t_+ +\O\big( \varepsilon^2 \big) \qquad \text{and} \qquad
\L(p_1,p_3)= \L(p_2,p_3)=0 \:. \]
Expanding the function~$\D(\vartheta)$ in~\eqref{DS3} in a quadratic Taylor polynomial,
we find that for a suitable real constant~$c$,
\beq
\begin{split}
\L(p_0,p_1)& = 1+c\, t_+^2 +\O \big( \varepsilon^3 \big) \\
\L(p_1,p_2)& = 1+c\, t_-^2 +\O \big( \varepsilon^3 \big) \\
\L(p_0,p_2)& = 1+c\, (t_+-t_-)^2 +\O\big( \varepsilon^3 \big) \:.
\end{split} \label{D-small}
\eeq

We choose a vector $u \in \C^4$ as
\begin{equation}\label{u-vector} 
u = \Big( 1, \;-1 + \frac{t_+}{t_-}, \;-\frac{t_+}{t_-},\; 0 \Big) 
\end{equation}
to find~$\langle u,Lu\rangle_{\C^4}=\O(\varepsilon^{3-2\beta})$, since
$|t_+/t_-|$ is of order $\lesssim \varepsilon^{-\beta}$.
Next, we choose the vector
\beq \label{vdef}
v=\alpha_1 \,u+\alpha_2 \,(0,0,0,1)
\eeq
for real numbers~$\alpha_1$ and~$\alpha_2$ to obtain 
\beq \label{Les}
\langle v,Lv\rangle_{\C^4}=\left\langle\left(\begin{array}{c} \alpha_1 \\
\alpha_2 \end{array}\right),\left(\begin{array}{cc}
\O(\varepsilon^{3-2\beta}) & b \,t_+ +\O(\varepsilon^2) \\[0.3em]
b \,t_+ +\O(\varepsilon^2) & 1 \end{array}\right)\left(\begin{array}{c} \alpha_1 \\
\alpha_2 \end{array}\right)\right\rangle.
\eeq
Since~$t_+^2 > \varepsilon^{2+2\beta}$ and~$\beta<1/6$,
the determinant of the matrix on the right is negative if~$\varepsilon$ is chosen sufficiently 
small, in contradiction to Corollary~\ref{corollary4}.
\QED
We remark that, for this argument, $\beta<\frac{1}{4}$ would suffice.

\subsection{Proof that~$q_1$ or~$q_2$ are Accumulation Points of~$M$}
\begin{Lemma} 
One of the points~$q_1$ or~$q_2$ is an accumulation point of~$M$.
\end{Lemma}
\Proof Assume conversely that~$q_1$ and~$q_2$ are not in the support of~$M$
or are isolated points in~$M$. We rotate the point~$p$ about the north pole by considering the curve
\[ p_s = (\cos s, \sin s, 0) \:. \]
Then there is~$s_0>0$ such that for all~$|s|< s_0$, the set~$\K(p_s) \cap M$ is empty.
Therefore, we can proceed just as in the proof of Lemma~\ref{lemma13} to
obtain a contradiction.
\QED

By symmetry, we may assume that~$q_1$ is an accumulation point.
We want to derive a contradiction.
According to Theorem~\ref{thm1}, the points in~$M$ in a neighborhood of~$q_1$
lie on a finite number of real analytic curves through~$q_1$.
We consider the following two cases: \label{caseab}
\begin{itemize}
\item[(a)] There is a real analytic curve~$c$ through~$q_1$ which osculates~$\K(p)$
and has the property that~$q_1$ is an accumulation point of support points on~$c$.
\item[(b)] There is no such curve.
\end{itemize}
Here we say that a curve {\em{osculates}} the light cone~$\K(p)$ if it approximates the
light cone to second order (i.e.\ if it is tangent to~$\K(p)$ in~$p$, and if its intrinsic curvature at~$p$
coincides with that of the geodesic circle~$\K(p)$).

We treat these two cases one after the other. We begin with the
simpler case~(b) (Section~\ref{secnoosculate}) and then consider case~(a)
(Section~\ref{secosculate}).

\subsection{Accumulation Point on Non-Osculating Curves through~$q_1$} \label{secnoosculate}
In this section, we shall rule out case~(b) above
(see Proposition~\ref{prpnob} below).
Thus we assume that every sequence of support points which converge
to~$q_1$ lies on a finite number of smooth curves, none of which osculates~$\K(p)$.
We consider any such curve and denote it by~$c$. Thus the function
\beq \label{nonosc}
c : (-\delta, \delta ) \rightarrow S^2
\eeq
is a smooth curve with the property that~$q_1$ is an accumulation point
of support points on~$c$. We parametrize by arc length and assume that~$c(0)=q_1$.
Moreover, this curve does not osculate~$\K(p)$.

Next, we let~$p_s$ be the curve
\beq \label{pscurve}
p_s = \begin{pmatrix} \cos s \\ \sin s \\ 0 \end{pmatrix}
\eeq
parametrizing the equator through~$p$.
Then the open light cone~$\I(p_s)$ around~$p_s$ may contain support points which are not
in~$\I(p)$. But, knowing that~$q_1$ and~$q_2$ are the only support points on~$\K(p)$,
the additional support points in~$\I(p_s)$ must lie in small balls centered at~$q_1$
and~$q_2$. In the next proposition, we quantify the radius~$r$ of these balls as a function of~$s$.
To this end, for~$r>0$ and~$q \in S^2$ we introduce the sets
\beq \label{Udeldef}
U_r = \big( B_r(q_1) \cup B_r(q_2) \big) \setminus \{q_1, q_2\} \qquad \text{and}
\qquad \I_r(q) = \I(q) \setminus U_r
\eeq
(where~$B_r$ denotes a geodesic ball in~$S^2$).
\begin{Lemma} \label{prpsize}
There are constants~$\kappa, s_0>0$ such that choosing
\beq \label{deltachoice}
r(s) = \kappa\, |s| \:,
\eeq
the equation
\beq \label{deltacond}
\I_{r(s_1)}(p) \cap M = \I_{r(s_1)}(p_s) \cap M
\qquad \text{holds for all~$s, s_1$ with~$|s| \leq s_1 \leq s_0$} \:.
\eeq
\end{Lemma}
\Proof
Since the center of the balls~$\I(p_s)$ depends continuously on~$s$, any point which lies
in~$\I(p_s)$ but not in~$\I(p)$ must lie in~$\K(p_{s'})$ for some~$s'$ with~$|s'| \leq s$.
The same is true for any point which lies in~$\I(p)$ but not in~$\I(p_s)$.
Therefore, it suffices to consider the boundary of the light cone~$\K(p_s)$ for small~$s$.

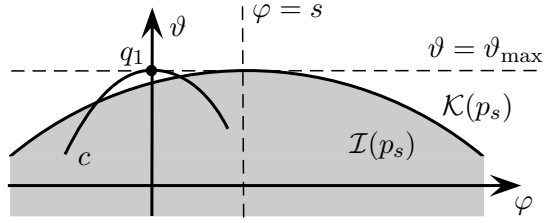
\begin{figure}
%
\psscalebox{1.0 1.0} 
{
\begin{pspicture}(-3,-1.4050229)(9.775,1.4050229)
\definecolor{colour0}{rgb}{0.8,0.8,0.8}
\psarc[linecolor=black, linewidth=0.04, fillstyle=solid,fillcolor=colour0, dimen=outer](3.175,-4.349977){4.9}{50.0}{130.0}
\psframe[linecolor=colour0, linewidth=0.02, fillstyle=solid,fillcolor=colour0, dimen=outer](6.33,-0.5899771)(0.03,-1.3899771)
\psline[linecolor=black, linewidth=0.04, arrowsize=0.093cm 5.0,arrowlength=1.44,arrowinset=0.4]{->}(0.0,-0.97997713)(6.825,-0.97997713)
\psline[linecolor=black, linewidth=0.04, arrowsize=0.093cm 5.0,arrowlength=1.44,arrowinset=0.4]{->}(1.915,-1.3849771)(1.92,1.3900229)
\rput[bl](2.15,1.0150229){\normalsize{$\vartheta$}}
\rput[bl](0.93,-0.6999771){\normalsize{$c$}}
\pscircle[linecolor=black, linewidth=0.02, fillstyle=solid,fillcolor=black, dimen=outer](1.91,0.5500229){0.075}
\psline[linecolor=black, linewidth=0.02, linestyle=dashed, dash=0.17638889cm 0.10583334cm](3.135,-1.4049771)(3.12,1.4050229)
\rput[bl](1.49,0.6250229){\normalsize{$q_1$}}
\psbezier[linecolor=black, linewidth=0.04](0.76,-0.5649771)(1.0841782,0.09335932)(1.4554658,0.5545028)(1.92,0.5550228881835938)(2.3845341,0.555543)(2.7357974,0.16600485)(2.925,-0.2249771)
\psline[linecolor=black, linewidth=0.02, linestyle=dashed, dash=0.17638889cm 0.10583334cm](0.045,0.5450229)(7.03,0.5500229)
\rput[bl](5.585,0.6700229){\normalsize{$\vartheta=\vartheta_{\max}$}}
\rput[bl](6.73,-1.3699771){\normalsize{$\varphi$}}
\rput[bl](3.25,1.195023){\normalsize{$\varphi=s$}}
\rput[bl](4.545,-0.6349771){\normalsize{$\I(p_s)$}}
\rput[bl](5.775,-0.1299771){\normalsize{$\K(p_s)$}}
\end{pspicture}
}
\caption{The geometry near~$q_1$ for a non-osculating curve.}
\label{figq}
\end{figure}

For any~$\varepsilon>0$, the set~$\K(p) \setminus (U_\varepsilon \cup \{q_1, q_2\})$
is compact and does not contain any points of~$M$ (see Proposition~\ref{lemmaq1q2}).
Thus it has an open neighborhood which does not intersect~$M$.
Consequently, there is~$s_0>0$ such that the sets~$\K(p_s) \setminus (U_\varepsilon \cup \{q_1, q_2\})$
do not intersect~$M$ for all~$s$ with~$|s| \leq s_0$. Therefore, it remains to show that
\[ 
\K(p_s) \cap \big( U_\varepsilon \setminus U_{r(s_1)} \big) \cap M = \varnothing \]
for all~$s,s_1$ as in~\eqref{deltacond}.
By symmetry, it clearly suffices to consider a neighborhood of the point~$q_1$ (see Figure~\ref{figq} in the case that the curve~$c$ is tangential to~$\K(p)$).

Parametrizing the sphere again by~$\vartheta \in (-\frac{\pi}{2}, \frac{\pi}{2})$ and~$\varphi \in (-\pi, \pi)$ with
\[ 
q(\vartheta, \varphi) = \begin{pmatrix} \cos \vartheta\: \cos \varphi \\  \cos \vartheta\: \sin \varphi \\ \sin \vartheta \end{pmatrix} \]
(thus~$\vartheta$ is the angle measured from the equator), the boundary of the light cone~$\K(p_s)$ is given by
\beq \label{thetaexp}
\vartheta = \vartheta_{\max} - a\, (\varphi-s)^2 + \O((\varphi-s)^4) \:.
\eeq
with~$a> 0$ (more precisely, $a=\cot(\vartheta_{\max})/2$; see the proof
of Lemma~\ref{lemma-osculate} below).

We first consider the case that the non-osculating curve~$c$ introduced in~\eqref{nonosc} intersects~$\K(p)$
transversally. In the subcase that the curve does not intersect~$\K(p)$ orthogonally, it can be parametrized by
\[ \vartheta = \vartheta_{\max} +b_1 \, \varphi + \O\big(\varphi^2\big) \]
with~$b_1 \neq 0$. Hence those points in~$M$ which lie in~$\K(p_s)$ are described
by the equation
\begin{align*}
\vartheta_{\max} - a\, (\varphi-s)^2 + \O((\varphi-s)^4) &= \vartheta_{\max} +b_1 \, \varphi + \O\big(\varphi^2\big) \\
\Longleftrightarrow \qquad
b_1 \, \varphi - 2 a s \,\varphi + as^2 &= \O \big( (\varphi-s)^4 \big) + \O \big(\varphi^2 \big) \:.
\end{align*}
For small~$s$, the implicit function theorem implies that this equation has solutions of the form
\beq \label{phies}
\varphi = -\frac{as^2}{b_1} + \O\big(s^3) \:.
\eeq
Using this expansion in~\eqref{thetaexp}, we also find that~$|\vartheta - \vartheta_{\max}|
\lesssim s^2$. Combining these estimates, we conclude that
for small~$s$, the points in~$M$ which are in~$\K(p_s)$
must lie inside a ball centered at~$q_1$ of radius~$r \lesssim s^2$.

In the subcase that the curve~$c$ intersects~$\K(p)$ orthogonally, for small~$|t|$
the curve~$c(t)$ is contained in the double cone
\[ \big| \vartheta - \vartheta_{\max} \big| \geq |\varphi| \:. \]
Therefore, the estimate~\eqref{phies} implies that the points on~$c$ which are in~$\K(p_s)$
satisfy the inequality
\[ \big|\varphi| \leq as^2 \:, \]
showing again that the points in~$M$ which are in~$\K(p_s)$
lie inside a ball centered at~$q_1$ of radius~$r \lesssim s^2$.

It remains to consider the case that the curve~$c$ does {\em{not}} intersect~$\K(p)$ transversally
or, in other words, that~$c$ is tangential to~$\K(p)$.
Then the non-osculating curve~$c$  can be parametrized by
\[ \vartheta = \vartheta_{\max} - \tilde{a}\, \varphi^2 + \O\big(\varphi^3\big) \]
with~$\tilde{a} \neq a$. Hence those points in~$M$ which lie on~$\K(p_s)$ are described
by the equations
\begin{align*}
\vartheta_{\max} - a\, (\varphi-s)^2 + \O((\varphi-s)^4) &= \vartheta_{\max} - \tilde{a}\, \varphi^2 + \O \big(\varphi^3 \big) \\
\Longleftrightarrow \qquad
(\tilde{a}-a)\, \varphi^2 + 2 a s \,\varphi - as^2 &= \O \big( (\varphi-s)^4 \big) + \O \big(\varphi^3 \big) \:.
\end{align*}
The resulting quadratic polynomial in~$\varphi$ has the roots
\[ \varphi_{1\!/\!2} = \frac{as}{\tilde{a}-a} \pm \frac{|s|}{\tilde{a}-a} \:\sqrt{2a^2 - a \tilde{a}
+ \frac{\tilde{a}-a}{s^2} \Big( \O \big( (\varphi-s)^4 \big) + \O \big(\varphi^3 \big) \Big) } \:. \]
This shows that for small~$s$, the points in~$M$ which are on~$\K(p_s)$
must lie inside a ball centered at~$q_1$ of radius~$r \lesssim |s|$.
\QED

In the next lemma we estimate derivatives of~$\D$ along the curve~$p_s$.
We begin with a preparatory lemma.
\begin{Lemma} Assume that~$\tau > \sqrt{6}$. Let~$p_s$ be the curve~\eqref{pscurve}.
Then for any~$q \in \I(p)$,
\begin{align}
\bigg|\frac{d \D(p_s, q)}{ds}\Big|_{s=0} \bigg| &\leq \frac 12\, 
\big(\tau^2+1 \big) \:\big| q_2 \big| \label{D1} \\
\frac{d^2 \D(p_s, q)}{ds^2}\Big|_{s=0} &< 0 \label{D2}
\end{align}
(where~$q_2$ denotes the second component of~$q \in \R^3$).
\end{Lemma}
\Proof We parametrize~$q$ as
\beq \label{qform}
q(\alpha, \varphi) = \begin{pmatrix} \sin \alpha\: \cos \varphi \\  -\sin \alpha\: \sin \varphi \\ \cos \alpha \end{pmatrix} \:.
\eeq
Since the angle between~$p_s$ and~$q$ depends only on~$\varphi+s$, we may set~$s$ to zero
and differentiate with respect to~$\varphi$.
By explicit computation, one finds that
\begin{align}
\cos \vartheta_{pq} &= \sin \alpha\: \cos \varphi \\
\D\big(p,q(\alpha, \varphi) \big) &= 
\frac 14\: (1+ \sin \alpha\: \cos \varphi) \left( 2 - \tau^2 \:(1 - \sin \alpha\: \cos \varphi) \right) \label{Dform} \\
\frac{\partial}{\partial \varphi} \D \big( p,q(\alpha, \varphi) \big)
&= -\frac 12 \,\big(1+\tau^2 \:\cos \varphi\: \sin \alpha\big) \:\sin \alpha\: \sin \varphi \:.\label{DDform}
\end{align}
Comparing the last equation with~\eqref{qform}, one immediately obtains the estimate~\eqref{D1}.

In order to prove the inequality~\eqref{D2}, we must show that the implication
\[ \D\big(p, q(\alpha, \varphi) \big)>0 \quad \Longrightarrow \quad
\frac{\partial^2}{\partial \varphi^2} \D \big( p,q(\alpha, \varphi) \big) <0 \]
holds for all~$\alpha \in \R$. From~\eqref{Dform}, one sees that
\[ 
\D \big(p,q(\alpha, \varphi) \big) > 0 \qquad \Longleftrightarrow \qquad \cos \varphi > \frac{\tau^2-2}{\tau^2\: \sin \alpha} 
\quad \text{and} \quad \sin \alpha\, \cos \varphi \neq -1 \:. \]
It follows that 
\begin{align*}
\frac{\partial^2}{\partial \varphi^2} \D \big( p,q(\alpha, \varphi) \big)
&= \frac 12 \sin \alpha \:\Big( -2 \cos^2 \varphi \: \tau ^2 \:\sin \alpha
-\cos \varphi +\tau ^2\: \sin \alpha \Big) \\
&\leq \frac{1}{2\tau^2} \:\big(-6+7 \tau^2-2 \tau^4+\tau^4 \:\sin^2 \alpha \big) 
\leq \frac{1}{2\tau^2} \:\big(-6+7 \tau^2-\tau^4 \big) \:.
\end{align*}
This is negative if~$\tau > \sqrt{6}$, concluding the proof.
\QED

After these preparations, we can now rule out case~(b) on page~\pageref{caseab}:
\begin{Prp} \label{prpnob}
If~$q_1$ is an accumulation point of~$M$, then it is an accumulation point of
points in~$M$ which lie on a curve which osculates~$\K(p)$ in~$q_1$.
\end{Prp}
\Proof We again proceed indirectly and assume that no such osculating curve exists.
Then the result of Lemma~\ref{prpsize} applies.
Choosing~$s_0$ and~$r(s)$ as in this lemma, by decomposing
the integral in~\eqref{ldef} we obtain
\[ \ell(p) = \int_{\I_{r(s_0)}(p)} \D(p, q) \: d\rho_{\min}(q) + 
\int_{U_{r(s_0)}} \L(p, q) \: d\rho_{\min}(q) \]
(where we used that~$\I_{r(s_0)}(p) \subset \I(p)$) and, similarly,
\begin{align*}
\ell(p_s) &= \int_{\I_{r(s_0)}(p_s)} \D(p_s, q) \: d\rho_{\min}(q) + 
\int_{U_{r(s_0)}} \L(p_s, q) \: d\rho_{\min}(q) \\
&= \int_{\I_{r(s_0)}(p)} \D(p_s, q) \: d\rho_{\min}(q) + 
\int_{U_{r(s_0)}} \L(p_s, q) \: d\rho_{\min}(q) \:,
\end{align*}
where in the last step we used~\eqref{deltacond} for~$s_1=s_0$ in order to change the integration
range from $\I_{r(s_0)}(p_s)$ to~$\I_{r(s_0)}(p)$. We thus obtain
\begin{align}
&\ell(p_s) - \ell(p) \nonumber \\
&= \int_{\I_{r(s_0)}(p)} \!\!\! \big( \D(p_s, q) - \D(p, q) \big) \: d\rho_{\min}(q)
+ \int_{U_{r(s_0)}} \!\!\! \big( \L(p_s, q) - \L(p, q) \big) \: d\rho_{\min}(q) \:. \label{elldecomp}
\end{align}

The last integral can be estimated as follows. For any~$q \in U_{r(s_0)}$
and sufficiently small~$s_0$,
\begin{align*}
\big| \L(p_s, q) - \L(p, q) \big| &\leq 
|s|\: \sup_{|s'|<|s|} \Big| \frac{d \D(p_{s'}, q)}{ds'} \Big|
= |s|\: \bigg( \bigg| \frac{d \D(p_{s}, q)}{ds} \Big|_{s=0} \bigg| + \O(s) \bigg) \\
&\!\!\!\overset{\eqref{D1}}{\leq} c(\tau)\: |s|\: r(s_0) + \O\big( s^2 \big)
\overset{\eqref{deltachoice}}{\leq} \tilde{c}(\tau)\: s^2\:.
\end{align*}
Integrating over~$U_{r(s_0)}$ gives the estimate
\[ \bigg| \int_{U_{r(s_0)}}  \big( \L(p_s, q) - \L(p, q) \big) \: d\rho_{\min}(q) \bigg|
\leq \tilde{c}(\tau)\: s^2\: \rho_{\min}\big(U_{r(s_0)} \big)\:. \]
Since in the limit~$s_0 \rightarrow 0$, the sets~$U_{r(s_0)}$ shrink to the empty set
(note that, according to our definition~\eqref{Udeldef}, the set~$U_{r(s_0)}$ do not contain the points~$q_1$
and~$q_2$),
for any~$\varepsilon>0$ we can arrange by decreasing~$s_0$ (and therefore~$r$) 
that~$\rho_{\min}(U_{r(s_0)})<\varepsilon$.
Since~$\varepsilon$ is arbitrary and~$|s| \leq s_0$, we conclude that
\[ 
\int_{U_{r(s_0)}}  \big( \L(p_s, q) - \L(p, q) \big) \: d\rho_{\min}(q) = o \big(s_0^2 \big) \:. \]
Clearly, this estimate applies just as well if in the integrand we replace~$\L$ by~$\D$.
This makes it possible to change the range of the first integral in~\eqref{elldecomp}
from~$\I_{r(s_0)}(p)$ to~$\I(p)$, i.e.
\[  \ell(p_s) - \ell(p) 
= \int_{\I(p)} \big( \D(p_s, q) - \D(p, q) \big) \: d\rho_{\min}(q)
+ o \big(s_0^2 \big) \:. \label{elloneint} \]

The obtained integral is obviously smooth in~$s$.
Performing a second order Taylor expansion in~$s$ gives
\begin{align*}
\ell \big(p_s \big) - \ell(p) = s\: \frac{d}{ds'} \dd_{\I(p)}\big(p_{s'}\big)\Big|_{s'=0}
+ \frac{s^2}{2} \:\frac{d^2}{ds'^2} \dd_{\I(p)}\big(p_{s'}\big)\Big|_{s'=0} + o\big(s_0^2 \big) \:,
\end{align*}
valid for all~$|s| \leq s_0$ and~$s_0$ sufficiently small,
where we used the notation
\[ \dd_U(x) := \int_U \D(x,y)\: d\rho_{\min}(y) \quad \in C^\infty(S^2)\:. 
\]
Choosing~$s=s_0$ and renaming~$s_0$ to~$s$, we obtain a Taylor expansion for~$\ell(p_s)$,
\[ \ell \big(p_s \big) = \ell(p) + s\: \frac{d}{ds'} \dd_{\I(p)}\big(p_{s'}\big)\Big|_{s'=0}
+ \frac{s^2}{2} \:\frac{d^2}{ds'^2} \dd_{\I(p)}\big(p_{s'}\big)\Big|_{s'=0} + o\big(s^2 \big) \:. \]
Now we can use the Euler-Lagrange equations~\eqref{EL1} to conclude that
\beq \label{ELcon}
\frac{d}{ds} \dd_{\I(p)}\big(p_{s}\big)\Big|_{s=0} = 0 \qquad \text{and} \qquad
\frac{d^2}{ds^2} \dd_{\I(p)}\big(p_{s}\big)\Big|_{s=0} \geq 0 \:.
\eeq

On the other hand, since~$\D$ is smooth and the integration range is compact, we know that
\begin{align*}
\frac{d^2}{ds^2}\: \dd_{\I(p)}\big(p_s) \Big|_{s=0} &=
\frac{d^2}{ds^2}\: \int_{\I(p)} \D(p_s, q)\: d\rho_{\min}(q) \Big|_{s=0} \\
&=\int_{\I(p)} \frac{d^2}{ds^2}\: \D(p_s, q)\: d\rho_{\min}(q) \Big|_{s=0} < 0\:,
\end{align*}
where in the last step we applied~\eqref{D2} together with the fact that~$\rho_{\min}$ is non-trivial on~$\I(p)$.
The last inequality contradicts~\eqref{ELcon}, concluding the proof.
\QED

\subsection{Accumulation Point on an Osculating Curve through~$q_1$} \label{secosculate}
In this section, we shall rule out case~(a) on page~\pageref{caseab} by proving
the following result:
\begin{Prp} \label{prpfinal}
Assume that~$p$ is a uniform two-sided accumulation point
on~$\gamma$ of scaling~$\beta$ with~$\beta$ according to~\eqref{betacond}.
Then there is no smooth curve through~$q_1$ which osculates~$\K(p)$
with the property that~$q_1$ is an accumulation point of support points on~$c$.
\end{Prp} \noindent
Proving this proposition and combining it with Proposition \ref{prphausdorffbeta}
will also conclude the proof of Theorem~\ref{thm0}.

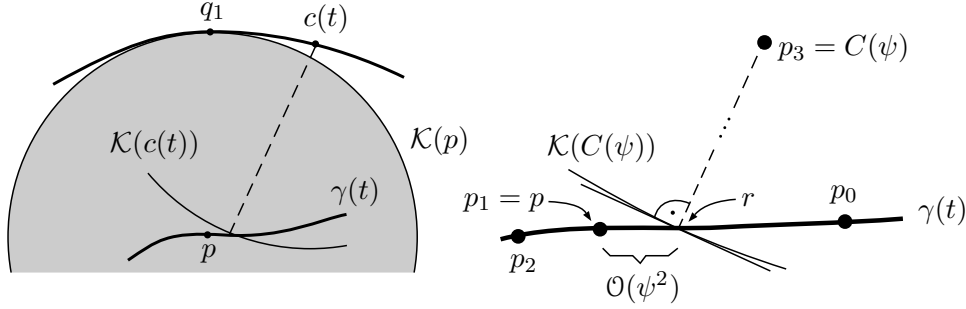
\begin{figure}
%
\psscalebox{1.0 1.0} 
{
\begin{pspicture}(0,-1.985)(12.82776,1.985)
\definecolor{colour0}{rgb}{0.8,0.8,0.8}
\psarc[linecolor=black, linewidth=0.02, fillstyle=solid,fillcolor=colour0, dimen=outer](2.8777592,-1.1116667){2.72}{-10.0}{190.0}
\rput[bl](2.7277594,-1.455){\normalsize{$p$}}
\rput[bl](2.7077594,1.775){\normalsize{$q_1$}}
\pscircle[linecolor=black, linewidth=0.06, fillstyle=solid,fillcolor=black, dimen=outer](8.028593,-1.0075){0.0975}
\pscircle[linecolor=black, linewidth=0.02, fillstyle=solid,fillcolor=black, dimen=outer](2.8452594,1.6175){0.0475}
\psbezier[linecolor=black, linewidth=0.04](0.75775933,0.92)(1.731664,1.4640204)(2.1978967,1.6095821)(2.8474996,1.62)(3.4971025,1.630418)(4.4074345,1.5068605)(5.4227595,1.0126786)
\psbezier[linecolor=black, linewidth=0.04](1.7527593,-1.42)(2.1415815,-1.1566635)(2.2670991,-1.0355201)(3.0416737,-1.085)(3.8162482,-1.1344799)(4.3166637,-0.874018)(4.657759,-0.81)
\rput[bl](4.447759,-0.705){\normalsize{$\gamma(t)$}}
\rput[bl](4.0877595,1.605){\normalsize{$c(t)$}}
\rput[bl](5.5027595,-0.02){\normalsize{$\K(p)$}}
\psarc[linecolor=black, linewidth=0.02, dimen=outer](4.174426,1.5483333){2.815}{-140.0}{-80.0}
\rput[bl](1.5127593,-0.075){\normalsize{$\K(c(t))$}}
\pscircle[linecolor=black, linewidth=0.02, fillstyle=solid,fillcolor=black, dimen=outer](4.2452593,1.4575){0.0475}
\psline[linecolor=black, linewidth=0.02, linestyle=dashed, dash=0.17638889cm 0.10583334cm](4.249426,1.4333333)(3.0960927,-1.0866667)
\psbezier[linecolor=black, linewidth=0.06](6.706016,-1.1411686)(7.6126404,-0.9028565)(8.840325,-1.0260257)(9.639446,-0.9821279186860966)(10.438567,-0.93823016)(11.235956,-0.9308823)(12.0545025,-0.8671647)
\psline[linecolor=black, linewidth=0.02, linestyle=dashed, dash=0.17638889cm 0.10583334cm](9.541093,0.051666666)(9.062759,-1.0133333)
\psline[linecolor=black, linewidth=0.02](10.292759,-1.5633334)(7.7727594,-0.41)
\psbezier[linecolor=black, linewidth=0.02](7.6077595,-0.21333334)(8.004541,-0.4463758)(8.628882,-0.80294675)(9.082389,-1.0028070175438757)(9.535895,-1.2026672)(10.06248,-1.3990703)(10.502759,-1.5366666)
\pscircle[linecolor=black, linewidth=0.06, fillstyle=solid,fillcolor=black, dimen=outer](6.9402595,-1.1025){0.0975}
\pscircle[linecolor=black, linewidth=0.06, fillstyle=solid,fillcolor=black, dimen=outer](11.285259,-0.9058333){0.0975}
\rput[bl](7.324426,-0.11){\normalsize{$\K(C(\psi))$}}
\rput[bl](6.217759,-0.77){\normalsize{$p_1=p$}}
\rput[bl](6.8327594,-1.585){\normalsize{$p_2$}}
\rput[bl](11.092759,-0.66){\normalsize{$p_0$}}
\rput[bl](12.257759,-0.955){\normalsize{$\gamma(t)$}}
\psline[linecolor=black, linewidth=0.02](8.062759,-1.26)(8.16276,-1.36)(8.462759,-1.36)(8.562759,-1.46)(8.637759,-1.36)(8.937759,-1.36)(9.062759,-1.26)
\rput[bl](8.094426,-1.985){\normalsize{$\O(\psi^2)$}}
\psbezier[linecolor=black, linewidth=0.02, arrowsize=0.05291667cm 2.0,arrowlength=1.4,arrowinset=0.0]{<-}(7.9127593,-0.9)(7.7727594,-0.6865306)(7.5627594,-0.64908165)(7.342759,-0.64)
\psbezier[linecolor=black, linewidth=0.02](9.214426,-0.66333336)(9.00917,-0.575307)(8.797456,-0.63929826)(8.751093,-0.8616666666666651)
\pscircle[linecolor=black, linewidth=0.02, fillstyle=solid,fillcolor=black, dimen=outer](2.8102593,-1.0825){0.0475}
\pscircle[linecolor=black, linewidth=0.02, fillstyle=solid,fillcolor=black, dimen=outer](8.995259,-0.7875){0.0275}
\psline[linecolor=black, linewidth=0.02, linestyle=dashed, dash=0.17638889cm 0.10583334cm](10.236093,1.5616666)(9.807759,0.6066667)
\psline[linecolor=black, linewidth=0.04, linestyle=dotted, dotsep=0.10583334cm](9.751093,0.47666666)(9.622759,0.21166667)
\pscircle[linecolor=black, linewidth=0.06, fillstyle=solid,fillcolor=black, dimen=outer](10.21526,1.4775){0.0975}
\rput[bl](10.40276,1.215){\normalsize{$p_3=C(\psi)$}}
\rput[bl](9.90776,-0.695){\normalsize{$r$}}
\psbezier[linecolor=black, linewidth=0.02, arrowsize=0.05291667cm 2.0,arrowlength=1.4,arrowinset=0.0]{<-}(9.20276,-0.91)(9.407496,-0.6786892)(9.619601,-0.6254788)(9.832759,-0.61)
\end{pspicture}
}
\caption{Accumulation point on an osculating circle.}
\label{figosculate}
\end{figure}

The proof is split up into several lemmas.
Clearly, we may assume that we are in case~(a) on page~\pageref{caseab}.
We again consider the smooth curve~$\gamma$
on which~$p$ is a uniform two-sided accumulation point of scaling~$\beta$
(see~\eqref{gammadef} and Definition~\ref{deftwosidedbeta}).
Moreover, we denote the osculating curve by~$c \in C^\infty((-\delta, \delta ),S^2)$. We again parametrize~$c$
by arc length and arrange that~$c(0)=q_1$ and~$\dot{c}(0)=(0,1,0)$ (see the left of Figure~\ref{figosculate}).
For our analysis, it is most convenient to parametrize these curves in the angle~$\varphi$:
\begin{Lemma} \label{lemmaphi}
There is~$\varphi_1>0$ and there are functions~$\Theta_C, \Theta_\Gamma \in C^\infty((-\varphi_1, \varphi_1), \R)$
such that the curves~$c$ and~$\gamma$ are given locally by
\beq \label{CG-phirep} 
C(\varphi)=\begin{pmatrix}
\cos\Theta_C(\varphi)\cos\varphi\\
\cos\Theta_C(\varphi)\sin\varphi\\
\sin\Theta_C(\varphi)\end{pmatrix}
\qquad \text{and} \qquad
\Gamma(\varphi)=\begin{pmatrix}
\cos\Theta_\Gamma(\varphi)\cos\varphi\\
\cos\Theta_\Gamma(\varphi)\sin\varphi\\
\sin\Theta_\Gamma(\varphi)\end{pmatrix} \:,
\eeq
respectively. Moreover,
\beq \label{inits}
\Theta_C(0) = \vartheta_{\max},\quad \Theta_\Gamma(0) = 0 \qquad \text{and} \qquad
\Theta_C'(0)= 0 = \Theta_\Gamma'(0) \:.
\eeq
\end{Lemma}
\Proof The method of proof is to construct the desired parametrization with the help of the
implicit function theorem.
We first represent the curves in spherical coordinates. Thus we choose angle functions
\begin{align*}
\vartheta_c,\vartheta_\gamma & : (-\delta,\delta)\to (-\pi/2,\pi/2) \\
\varphi_c,\varphi_\gamma & : (-\delta,\delta)\to (-\pi,\pi) \:,
\end{align*}
such that 
\begin{equation}\label{rep-spherical}
c(t)=\begin{pmatrix} \cos\vartheta_c(t)\cos\varphi_c(t)\\
\cos\vartheta_c(t)\sin\varphi_c(t)\\
\sin\vartheta_c(t)\end{pmatrix} \quad\textnormal{and}\quad
\gamma(t)=\begin{pmatrix} \cos\vartheta_\gamma(t)\cos\varphi_\gamma(t)\\
\cos\vartheta_\gamma(t)\sin\varphi_\gamma(t)\\
\sin\vartheta_\gamma(t)\end{pmatrix}
\end{equation}
for all~$t\in (-\delta,\delta)$.
In particular, the third components $c^3(\cdot)$ and~$\gamma^3(\cdot)$
of these curves are smooth.
Moreover, evaluating at~$t=0$, we
find that~$\vartheta_c(0)=\vartheta_{\max}$ and $\vartheta_\gamma(0)=0$
and thus~$c^3(0), \gamma^3(0)\in (-1,1)$.
By choosing~$\delta$ sufficiently small, we can arrange that~$c^3(t), \gamma^3(t)\in (-1,1)$
for all~$t \in (-\delta, \delta)$.
Hence we may apply the $\arcsin$ to obtain
functions
\[ \vartheta_c, \vartheta_\gamma \in C^\infty \big( (-\delta, \delta) \big) \]
and
\[ \dot{\vartheta}_c(0) = \frac{\dot{c}^3(0)}{\sqrt{1-(c^3(0))^2}} =0\:,\qquad
\dot{\vartheta}_\gamma(0) = \frac{\dot{\gamma}^3(0)}{\sqrt{1-(\gamma^3(0))^2}} =0 \:. \]
Choosing $\delta>0$ even smaller if necessary,
we can use that~$\vartheta_c(0)=\vartheta_{\max}$ and $\vartheta_\gamma(0)=0$ to arrange that both
functions $\cos\vartheta_c(t)$ and $ \cos\vartheta_\gamma(t)$ are strictly positive for all $|t|\le \delta$.
As a consequence, we can infer from the second components $c^2$ and $\gamma^2$ in 
\eqref{rep-spherical} that 
$$
\sin\varphi_c(\cdot)=\frac{c^2(\cdot)}{\cos\vartheta_c(\cdot)}\in
C^\infty((-\delta,\delta)) \quad\text{and} \quad \varphi_c(0)=0\:,
$$
and analogously, $\sin\varphi_\gamma(\cdot)$ is smooth
with $\varphi_\gamma(0)=0$.
If necessary, we choose $\delta$ even smaller so that
$\sin\varphi_c(\cdot)$ and $\sin\varphi_\gamma(\cdot)$ take values
in $(-1,1)$, which is possible due to continuity and the initial
values $c^2(0)=0=\gamma^2(0)$. We conclude that
\[ \varphi_c, \varphi_\gamma \in C^\infty((-\delta,\delta)) \:. \]

Now we use our assumptions on the tangent vectors~$\dot{c}(0)$ and $\dot{
\gamma}(0)$ to derive from \eqref{rep-spherical} by means of the
initial values $\vartheta_c(0)=\vartheta_{\max}$ and $\varphi_c(0)=0$ that
\begin{align}
1=\dot{c}^2(0) & \overset{\eqref{rep-spherical}}{=}
-\sin\vartheta_c(0)\cdot\dot{\vartheta}_c(0) \: \sin\varphi_c(0)+
\cos\vartheta_c(0) \:\cos\varphi_c(0)\cdot\dot{\varphi}_c(0) \notag \\
&\;\,= \cos \vartheta_{\max} \cdot\dot{\varphi}_c(0) \:. \label{dotcphi}
\end{align}
This implies that~$\dot{\varphi}_c(0)=1/\cos \vartheta_{\max} >0$.
The inverse function theorem yields a strictly monotone
function~$t_c \in C^\infty(-\varphi_1,\varphi_1))$ with~$\varphi_1>0$ such
that~$t_\pm := t_c(\pm \varphi_1) \in (-\delta, \delta)$ and
\begin{align*}
t_c(\varphi_c(t)) &= t \qquad\; \text{for all~$t\in (t_-,t_+)$} \\
\varphi_c(t_c(\varphi)) &= \varphi \qquad \text{for all~$\varphi\in (-\varphi_1,\varphi_1)$} \\
\frac{d}{d\varphi} t_c(\varphi) \Big|_{\varphi=0}&= t_c'(0)=\frac{1}{\dot{\varphi_c}(0)}=\cos \vartheta_{\max}.
\end{align*}
Setting $C(\varphi):=c\circ t_c(\varphi)$ and $\Theta_C(\varphi):=
\vartheta_c\circ t_c(\varphi)$ for $\varphi\in (-\varphi_1,\varphi_1)$
gives the representation of the curve~$c$ by a function~$C$ as in~\eqref{CG-phirep}.
In particular, the function $\Theta_C(\cdot)$ is smooth
and has the properties~$\Theta_C(0)=\vartheta_{\max}$ and
\[ 
\Theta_C'(0)=\left.\frac{d}{d\varphi}\right|_{\varphi=0}\Theta_C(\varphi)=\dot{\vartheta}_c(0)\cdot t'_c(0) =0 \:. \]

Analogously, one obtains a representation of the curve~$\gamma$ 
by a function~$\Gamma$ as in~\eqref{CG-phirep}
depending on $\varphi $ on the same parameter interval (if $\varphi_1$ is chosen sufficiently small).
Since~$\gamma(t)$ is parametrized
by arc length and in~$p$ is tangential to the equator, in analogy to~\eqref{dotcphi}
we obtain
\beq \label{phiGdot}
\dot{\varphi}_\gamma(0)=1 \:.
\eeq
Also here, $\Theta_\Gamma(\cdot)$ is smooth, this time with the initial value
$\Theta_\Gamma(0)=0$ and with $\Theta_\Gamma'(0)=0$. 
\QED

We now specify what ``osculating'' means in our new parametrization:
\begin{Lemma} \label{lemma-osculate}
In the parametrization~\eqref{CG-phirep}, the curve~$c$ osculates~$\K(p)$
if and only if the function~$\Theta_C(\varphi)$ satisfies the relation
\begin{equation}\label{osculating}
\Theta_C''(0)=-\cot \vartheta_{\max} \:.
\end{equation}
\end{Lemma}
\Proof In order to verify~\eqref{osculating}, we consider the
parametrization of $\mathcal{K}(p)$ by arc length
$$
k(t):=\begin{pmatrix}
\sqrt{1-\cos^2 \vartheta_{\max}} \\
\sin \vartheta_{\max}\cdot\sin(t/\sin \vartheta_{\max}) \\
\sin \vartheta_{\max}\cdot\cos(t/\sin \vartheta_{\max})\end{pmatrix}
$$
with $k(0)=q_1$ and $\dot{k}(0)=(0,1,0)$.
Using spherical coordinates $\vartheta_k(t),\varphi_k(t)$ as in~\eqref{rep-spherical} to
describe $k$, and writing the third component of~$k(t)$ as~$k^3(t)=\sin\vartheta_k(t)$, we obtain
$$
\ddot{\vartheta}_k(0)=-\frac{1}{\sin \vartheta_{\max}\cdot\cos \vartheta_{\max}}.
$$
Rewriting also $\vartheta_k$ as a function of $\varphi$ (similarly as done
for the general curve $c$ above), we obtain
$\Theta_K(\varphi):=(\vartheta_k\circ t_k)(\varphi)$ and thus
$$
\Theta_K''(\varphi)= \big( \ddot{\vartheta}_k\circ t_k \big)(\varphi) \:(t_k'(\varphi))^2+
\big(\dot{\vartheta}_k\circ t_k \big)(\varphi)\:t_k''(\varphi) \:.
$$
It follows that
$$
\Theta_K''(0)=\ddot{\vartheta}_k(0)\cos^2 \vartheta_{\max}=-\cot \vartheta_{\max}=\Theta_C''(0) \:,
$$
in agreement with~\eqref{osculating}.
\QED

We now consider the curves~$C(\psi)$ and~$\Gamma(\phi)$
for two different angles~$\psi, \phi \in (-\varphi_1, \varphi_1)$.
The following lemma shows that, for sufficiently small~$\psi$,
the curve~$\Gamma(\phi)$ intersects the light cone around~$C(\psi)$ transversally.
\begin{Lemma} \label{lem:transversal-intersection}
There are~$\psi_2, \phi_2 \in (0, \varphi_1)$ such that the following statements hold:
For every~$\psi \in (-\psi_2, \psi_2)$
there is~$\phi \in (-\phi_2, \phi_2)$ such
such that~$\Gamma(\phi) \in \K(C(\psi))$.
Moreover, the curve~$\Gamma$ intersects the light cone~$\K(C(\psi))$
transversally (see the right of Figure~\ref{figosculate}).
The angle~$\phi$ satisfies the scaling
\[ 
\phi = \O\big( \psi^2 \big) \:. \]
\end{Lemma}
\Proof 
We denote the angle between the points~$C(\psi)$ and~$\Gamma(\phi)$
by~$\vartheta_{\psi, \phi}$. Then, taking the scalar product between
the vectors in~\eqref{CG-phirep} and applying the sum rules, we obtain
\beq \label{cosformel}
\cos \vartheta_{\psi, \phi} = \cos(\phi-\psi)\: \cos\big(\Theta_C(\psi) \big)\:
\cos \big( \Theta_\Gamma(\phi) \big) + \sin\big(\Theta_C(\psi) \big)\:
\sin \big( \Theta_\Gamma(\phi) \big) \:.
\eeq
Expanding in a Taylor polynomial and using~\eqref{inits} and~\eqref{osculating},
we obtain
\beq \begin{split}
\cos \vartheta_{\psi, \phi} & = \cos \vartheta_{\max} 
+ \psi \,\phi\, \cos \vartheta_{\max}
-\frac{\psi^3}{6} \:\Theta_C^{(3)}(0)\:\sin \vartheta_{\max} \\
&\quad\: + \O\big( \phi^2 \big) + \O\big( \phi\, \psi^3 \big)+ \O\big( \psi^4 \big)  \:. \label{phiexp}
\end{split}
\eeq
The condition~$\Gamma(\phi) \in \K(C(\psi))$ is satisfied if and only
if~$\cos \vartheta_{\psi, \phi} = \cos \vartheta_{\max}$. Therefore, the expansion~\eqref{phiexp}
suggests a solution~$\phi$ of the form
\[ \phi  = \frac{\psi^2}{6} \:\Theta_C^{(3)}(0)\:\tan \vartheta_{\max} 
+ \O\big( \psi^3 \big) \:. \]
However, since the error terms in~\eqref{phiexp} also involve the unknown~$\phi$, the
existence of this solution is not obvious. Therefore, we prove existence using the following
fixed-point argument: In order to relate solutions of the equation~$\cos \vartheta_{\psi, \phi} = \cos \vartheta_{\max}$
to fixed points of a mapping~$T$, given~$\psi \neq 0$ and choosing
\[ \phi_2 = \frac{\psi^2}{3} \:\Theta_C^{(3)}(0)\:\tan \vartheta_{\max}\:, \]
we introduce the function~$T(\phi)$ by
\[ T \::\: [-\phi_2, \phi_2] \rightarrow \R \:,\qquad
T(\phi) = \phi - \frac{1}{\psi\: \cos \vartheta_{\max}}\:\big( \cos \vartheta_{\psi, \phi} - \cos \vartheta_{\max} \big) \:. \]
A straightforward computation using~\eqref{phiexp} yields that this function and its derivative have the Taylor expansions
\[ T(\alpha \phi_2) = \frac{\phi_2}{2}\: + \O (\psi) \qquad \text{and} \qquad
T'(\alpha \phi_2) = \O(\psi) \:, \]
uniformly in~$\alpha \in [-1,1]$. It follows that, for sufficiently small~$\psi$,
the function~$T$ maps the interval~$[-\phi_2, \phi_2]$ into itself.
Moreover, the mean value inequality gives the estimate
\[ \big|T(\phi) - T(\phi') \big| \leq |\phi-\phi'|\: \sup_{\alpha \in [-1,1]} \big| T'(\alpha \phi_2) \big|
= |\phi-\phi'|\: \O(\psi) \]
proving that, again for sufficiently small~$\psi$, the mapping~$T$ is a contraction.
Therefore, the Banach fixed-point theorem applies and yields the
desired~$\phi \in [-\phi_2, \phi_2]$ with~$\phi=T(\phi)$.
\QED

We now invoke our final argument. Since~$q_1$ is an accumulation point
of support points on the curve~$c$, we can choose~$\psi \in (-\psi_2, \psi_2)$
with~$|\psi|$ arbitrarily small such that~$C(\psi) \in M$.
We choose
\beq \label{epschoice}
\varepsilon = \big( \mathfrak{c} \,\psi^2 \big)^{\frac{1}{1+2\beta}}
\eeq
with a positive constant~$\mathfrak{c}$ to be determined below.
Note that by decreasing~$|\psi|$ we can make~$\varepsilon$
arbitrarily small. Since~$p$ is a uniform two-sided accumulation point on~$\gamma$
with scaling~$\beta$ (see Definition~\ref{deftwosidedbeta}), 
we can arrange by decreasing~$|\psi|$ that
there are parameter values~$t_\pm$ in the range~\eqref{uniformbeta}
such that the points~$\gamma(t_\pm)$ are in~$M$, i.e.
\[ \big( \mathfrak{c} \,\psi^2\big)^\frac{1+\beta}{1+2 \beta}
< |t_\pm| < \big( \mathfrak{c} \,\psi^2 \big)^{\frac{1}{1+2\beta}} \:. \]

We now improve the method used in the proof of Proposition~\ref{lemmaq1q2}.
Namely, we again consider the Gram matrix~$L$, but now for the four points
\[ p_0 =\gamma(t_+)\:,\qquad p_1 =p \:,\qquad p_2 =\gamma(t_-)
\qquad \text{and} \qquad p_3 =C(\psi) \:. \]
By changing the orientation of~$\gamma$ and decreasing~$|\psi|$ if necessary,
we can arrange that~$p_0$ lies inside and~$p_2$ lies outside
the light cone around~$C(\psi)$ (as shown on the right of Figure~\ref{figosculate};
the point~$p_1$, however, could lie inside or outside).
We denote the point where the curve~$\gamma$ intersects the boundary of this light cone
by~$r=\gamma(t_1)$ (see Figure~\ref{figosculate}; note that this point does in general not lie in~$M$).

In the previous estimates, we often parametrized the curve~$\gamma$
by angle~$\phi$. Since~$\gamma(t)$ is parametrized
by arc length and in~$p$ is tangential to the equator (see~\eqref{phiGdot}), 
we can bound~$\phi$ from above and below by~$t$,
implying that all the scalings for~$\phi$ hold just as well for~$t$.
In particular, Lemma~\ref{lem:transversal-intersection} implies that
\[ t_1 = \O(\psi^2) \:. \]
We next compute the first derivative of~$\D$ on the light cone in the direction of~$\gamma$.
Again parametrizing~$\gamma$ in the variable~$\phi$ and setting~$\phi_1 = \varphi_\gamma(t_1)$,
we obtain
\[ \frac{\partial}{\partial \phi} 
\D \big( \vartheta_{\psi, \phi} \big)\Big|_{\phi_1}
= -\frac{\D'\big( \vartheta_{\max} \big)}{\sin \vartheta_{\max}}\:
\frac{\partial \cos \vartheta_{\psi, \phi}}{\partial \phi} \Big|_{\phi_1} 
= -\frac{\D'\big( \vartheta_{\max} \big)}{\sin \vartheta_{\max}}\;
\psi \, \cos \vartheta_{\max} + \O\big(\psi^2\big) \:, \]
where in the last step we differentiated~\eqref{cosformel} and used~\eqref{inits} and~\eqref{osculating}.
Using again~\eqref{phiGdot}, we conclude that
\[ \frac{d}{dt} \D \big( \gamma(t), p_3 \big)\Big|_{t_1}
= g\, \psi + \O(\psi^2) \qquad \text{with $g \neq 0$} \:. \]

The last formula gives rise to the following expansions of the function~$\D$,
\begin{align*}
\D(p_0, p_3) &= (t_+-t_1) \: \frac{d}{dt}\D \big( \gamma(t),p_3 \big)\Big|_{t_1} 
+ \O \big( (t_+-t_1)^2 \big) \\
&= |g\,\psi|\: \big( \mathfrak{c} \,\psi^2\big)^\frac{1+\beta}{1+2 \beta}
\big(1 + \O(\psi)\big) \\
\D(p_1, p_3) &= -t_1 \: \frac{d}{dt}\D \big( \gamma(t),p_3 \big)\Big|_{t_1} 
+ \O \big( t_1^2 \big) = g \,\O \big(\psi^3 \big) \:.
\end{align*}
Using again the expansion~\eqref{D-small} with the same vectors~$u$
and~$v$ as in~\eqref{u-vector} and~\eqref{vdef},
we find $\langle u,Lu\rangle = O(\varepsilon^{3-2\beta})$, so that~\eqref{Les} becomes
\beq \label{Les1}
\langle v,Lv\rangle_{\C^4}=\left\langle\left(\begin{array}{c} \alpha_1 \\
\alpha_2 \end{array}\right),\left(\begin{array}{cc}
\O(\varepsilon^{3-2\beta}) & P(\psi,t_+,t_-)  \\[0.3em]
P(\psi,t_+,t_-) & 1 \end{array}\right)\left(\begin{array}{c} \alpha_1 \\
\alpha_2 \end{array}\right)\right\rangle,
\eeq
where the off-diagonal entries scale like
\beq \label{Ptdef}
P\big(\psi,t_+,t_- \big) = 
|g| \,\mathfrak{c}^\frac{1+\beta}{1+2 \beta}
 \,|\psi|\: \big(  \,\psi^2\big)^\frac{1+\beta}{1+2 \beta}\,
\big(1 + \O(\psi)\big)  +  g \:\frac{t_+-t_-}{t_-}\:\O(\psi^3) \:.
\eeq
According to~\eqref{uniformbeta}, the quotient~$t_+/t_-$ is bounded by
\[ \Big| \frac{t_+}{t_-} \Big| \leq \varepsilon^{-\beta} = 
\big( \mathfrak{c} \,\psi^2\big)^{-\frac{\beta}{1+2 \beta}} \:. \]
Using this bound in~\eqref{Ptdef} gives
\[ P\big(\psi,t_+,t_- \big)  = 
|g| \, \mathfrak{c}^\frac{1+\beta}{1+2 \beta} \: |\psi|^\frac{3+4\beta}{1+2 \beta}
\big(1 + \O(\psi)\big)  +  g \:\mathfrak{c}^{-\frac{\beta}{1+2 \beta}} \,
\O \Big( |\psi|^{\frac{3 + 4 \beta}{1+2 \beta}} \Big) \:. \]
Therefore, by increasing~$\mathfrak{c}$ we can arrange that the first term dominates,
so that for small~$|\psi|$,
\[ \big| P\big(\psi,t_+,t_- \big) \big| \geq
\frac{|g|}{2} \, \mathfrak{c}^\frac{1+\beta}{1+2 \beta} \: |\psi|^\frac{3+4\beta}{1+2 \beta} \:. \]

The previous estimates are not good enough for concluding that
the determinant of the matrix in~\eqref{Les1} is negative (which in turn is needed
in order to obtain a contradiction to Corollary \ref{corollary4}).
Namely, the product of the diagonal terms is of the order
\beq
\O\big( \varepsilon^{3-2\beta} \big) \overset{\eqref{epschoice}}{=}
\O\Big( \psi^{\frac{6-4\beta}{1+2 \beta}} \Big)  \:, \label{diagonal}
\eeq
whereas the product of the off-diagonal terms scales like
\beq \label{offdiagonal}
\big|P\big(\psi,t_+,t_- \big) \big|^2 \gtrsim |\psi|^\frac{6+8\beta}{1+2 \beta} \:.
\eeq
This explains why we must improve the error term in the upper left
matrix entry as follows:
\begin{Lemma} \label{lemmaimproved}
The error terms in the expansion~\eqref{D-small} can be improved to
\begin{align*}
\L(p_0,p_1)& = 1+c\, t_+^2 +\O \big( \varepsilon^4 \big) \\
\L(p_1,p_2)& = 1+c\, t_-^2 +\O \big( \varepsilon^4 \big) \\
\L(p_0,p_2)& = 1+c\, (t_+-t_-)^2 +\O\big( \varepsilon^4 \big) \:.
\end{align*}
\end{Lemma}
\Proof Intuitively speaking, the error terms come about for two reasons. First,
due to the higher derivatives of the function~$\D(\vartheta)$. But since~$\D$ is
even in~$\vartheta$, a Taylor expansion about~$\vartheta=0$ only involves
even derivatives. Therefore, the resulting error is of the order~$\varepsilon^4$.
The second reason for the error terms is the fact that the curve~$\gamma$
may have an intrinsic curvature. However, this error is even in the parameter~$t$
of the curve~$\gamma$, explaining the error term~$\O(\varepsilon^4)$.

In order to compute the effect of the intrinsic curvature in detail, we parametrize the
curve~$\gamma$ as (see again Figure~\ref{figsphere})
\beq \label{gamma3}
\gamma(t) = \big( \sqrt{1-y(t)^2-z(t)^2}, \;y(t), \;z(t) \big) \:.
\eeq
By a rotation around the $x$-axis, we can arrange that
\beq \label{init1}
y(0)=z(0)= \dot{z}(0) = 0 \:.
\eeq
The fact that the curve~$\gamma$ is parametrized by arc length implies that
\beq \label{init2}
\dot{y}(0)=1 \qquad \text{and} \qquad \ddot{y}(0) = 0 \:.
\eeq
Denoting the angle between the points~$\gamma(t)$ and~$\gamma(t')$
by~$\vartheta_{t,t'}$, the cosine of this angle can be computed conveniently
by taking the scalar products of the vectors in~$\R^3$ in the parametrization~\eqref{gamma3}.
A straightforward computation using~\eqref{init1} and~\eqref{init2} yields
\begin{align*}
\cos \vartheta_{0,t_+} &= 1- \frac{t_+^2}{2} + \O\big( \varepsilon^4 \big) \:,\qquad
\cos \vartheta_{t_-,0} = 1+\frac{t_-^2}{2} + \O\big( \varepsilon^4 \big) \\
\cos \vartheta_{t_-,t_+} &= 1+\frac{(t_+-t_-)^2}{2} + \O\big( \varepsilon^4 \big) \:.
\end{align*}
Since~$\D$ is a polynomial in~$\cos \vartheta$ (see~\eqref{DS3}),
we obtain the desired expansion with error term~$\O(\varepsilon^4)$.
\QED

\Proof[Proof of Proposition~\ref{prpfinal}]
Using the improved estimate of Lemma~\ref{lemmaimproved}, 
the expectation value in~\eqref{Les1} is replaced by
\[  \langle v,Lv\rangle_{\C^4}=\left\langle\left(\begin{array}{c} \alpha_1 \\
\alpha_2 \end{array}\right),\left(\begin{array}{cc}
\O(\varepsilon^{4-2\beta}) & P(\psi,t_+,t_-) \\[0.3em]
P(\psi,t_+,t_-) & 1 \end{array}\right)\left(\begin{array}{c} \alpha_1 \\
\alpha_2 \end{array}\right)\right\rangle \:. \]
Computing the
determinant of the matrix, the product off-diagonal terms scale again according to~\eqref{offdiagonal}.
However, the scaling of the product of the diagonal terms is improved from~\eqref{diagonal} to
\[ \O\big( \varepsilon^{4-2\beta} \big) =
\O\Big( \psi^{\frac{8-4\beta}{1+2 \beta}} \Big)  \:. \]
Therefore, using the assumption~$\beta<1/6$, the determinant is negative for small~$|\psi|$.
We thus obtain the desired contradiction to Corollary \ref{corollary4}.
\QED

\Thanks {{\em{Acknowledgments:}}
We would like to thank Tobias Kaiser for helpful discussions.
We are grateful to the referee for helpful comments on the manuscript.
F.F.\ would like to thank the Max Planck Institute for Mathematics in the Sciences in Leipzig for
hospitality while working on the manuscript. H.~vdM.'s work is partially funded by the
Excellence Initiative of the German federal and state governments.

\providecommand{\bysame}{\leavevmode\hbox to3em{\hrulefill}\thinspace}
\providecommand{\MR}{\relax\ifhmode\unskip\space\fi MR }
\providecommand{\MRhref}[2]{%
  \href{http://www.ams.org/mathscinet-getitem?mr=#1}{#2}
}
\providecommand{\href}[2]{#2}

\end{document}